\documentclass[10pt,twocolumn,twoside,journal]{IEEEtran}

\usepackage{graphicx} 
\usepackage{epstopdf}
\usepackage{amsmath} 
\usepackage{amssymb}  
\usepackage{amsthm}
\usepackage{xcolor}

\newcommand{\bbR}{\ensuremath{\mathbb{R}}}

\newcommand{\bigO}[1]{\ensuremath{\mathcal{O}\left( #1 \right)}}
\newcommand{\mcN}{\ensuremath{\mathcal{N}}}

\newcommand{\bfH}{\ensuremath{\mathbf{H}}}

\newcommand{\bfJ}{\ensuremath{\mathbf{J}}}
\newcommand{\bfm}{\ensuremath{\mathbf{m}}}
\newcommand{\bfn}{\ensuremath{\mathbf{n}}}
\newcommand{\bfr}{\ensuremath{\mathbf{r}}}

\newcommand{\bfx}{\ensuremath{\mathbf{x}}}
\newcommand{\bfX}{\ensuremath{\mathbf{X}}}
\newcommand{\bfy}{\ensuremath{\mathbf{y}}}

\newcommand{\bftheta}{\ensuremath{\boldsymbol \theta}}

\newcommand{\bfmu}{\ensuremath{\boldsymbol \mu}}

\newcommand{\muBar}{\ensuremath{\bar{\mu}_0}}
\newcommand{\deltaBar}{\ensuremath{\bar{\delta}^2_0}}
\newcommand{\deltaMu}{\ensuremath{\Delta_{\mu}}}
\newcommand{\deltaDelta}{\ensuremath{\Delta_{\delta}}}
\newcommand{\oneVec}[1]{\ensuremath{\mathbf{1}_{ #1 }}}

\DeclareMathOperator*{\Expectation}{\mathbb{E}}

\newcommand{\map}[3]{#1: #2 \rightarrow #3}
\newcommand{\trans}[1]{\left. #1 \right.^T}
\newcommand{\indicator}[1]{\ensuremath{\mathbf{1}\! \left( #1 \right)}}
\renewcommand{\Pr}[1]{\ensuremath{\mathrm{Pr} \left[ #1 \right]}}
\newcommand{\E}[1]{\ensuremath{\Expectation\left[ #1 \right]}}

\DeclareMathOperator*{\argmax}{arg\,max}

\DeclareMathOperator{\sumOp}{sum}

\DeclareMathOperator{\diag}{diag}

\def \etal {\emph{et al.}}

\newcommand{\beq}{\begin{equation}}
\newcommand{\eeq}{\end{equation}}
\newcommand{\bal}{\begin{align}}
\newcommand{\eal}{\end{align}}

\newtheorem{lem}{Lemma}
\newtheorem{thm}{Theorem}

\theoremstyle{definition}
\newtheorem{defn}{Definition}
\newtheorem{example}{Example}
\newtheorem{prob}{Problem}

\title{\LARGE \bf
Parameter estimation in softmax decision-making models with linear objective functions
}

\author{Paul Reverdy$^1$ and Naomi E. Leonard$^2$
\thanks{This research has been supported in part by ONR grants N00014-09-1-1074 and N00014-14-1-0635, and was conducted with Government support under and awarded by DoD, Air Force Office of Scientific Research, National Defense Science and Engineering Graduate (NDSEG) Fellowship, 32 CFR 168a.}
\thanks{$^1$Department of Electrical and Systems Engineering, University of Pennsylvania, Philadelphia, PA 19104, USA ({\tt\small{preverdy@seas.upenn.edu}}),

$^2$Department of Mechanical and Aerospace Engineering, Princeton University, Princeton, NJ 08544, USA ({\tt\small{naomi@princeton.edu}})}%
}

\begin{document}

\maketitle
\thispagestyle{empty}
\pagestyle{empty}

\begin{abstract}

With an eye towards human-centered automation, we contribute to the development of a systematic means to infer features of human decision-making from behavioral data.  Motivated by the common use of softmax selection in models of human decision-making, we study the maximum likelihood parameter estimation problem for softmax decision-making models with linear objective functions.  We present conditions under which the likelihood function is convex.  These allow us to provide sufficient conditions for convergence of the resulting maximum likelihood estimator and to construct its asymptotic distribution. In the case of models with nonlinear objective functions, we show how the estimator can be applied by linearizing about a nominal parameter value. We apply the estimator to fit the stochastic UCL  (Upper Credible Limit) model of human decision-making to human subject data. We show statistically significant differences in behavior across related, but distinct, tasks.

\end{abstract}

Note to Practitioners:
\begin{abstract} We propose and demonstrate a rigorous method to estimate parameters of softmax decision-making models.  These decision-making models hold great promise for use in developing model-based human-centered automation. We are motivated by the recently derived UCL (Upper Credible Limit) model, which predicts the choices that humans are likely to make when deciding among  alternatives with uncertain rewards. Key parameters of the model represent the human's intuition about the task, and estimating these parameters from behavioral data would allow an automated system to learn about its human supervisor.  Our parameter estimation method is fast enough to be implemented in real time for most scenarios, although our analysis of the method holds when the model has a particular linear structure. We show how to extend the method to a more general nonlinear model using linearization, and we show that the linearization approach works for the motivating UCL model. The parameter estimation method with linearization can be used for other nonlinear models; however, the domain of its validity may vary.
\end{abstract}

Primary and Secondary Keywords
\begin{IEEEkeywords}
Primary Topics: Estimation, Automation, Decision-Making
\end{IEEEkeywords}


\section{Introduction}

In a variety of decision-making scenarios an agent selects one among a discrete set of options $i \in \{1, \ldots, m \}$ and receives a reward associated with the selection.   The agent's goal is to make a selection or a sequence of selections to maximize reward.  For example, a human air traffic controller selects among options for allocating aircraft for takeoff, and the reward is a measure of efficiency of flight departures associated with the selected option \cite{VR-HB:11}. Often the decision-making task is challenging, especially when there is uncertainty or there are complex dependencies associated with options and rewards, as in the air traffic control example.   In this paper we propose a rigorous method to estimate features of humans decision-making that can be used to enable human-centered automation.

Much research has gone into studying how humans decide among options and what conditions lead to good decision-making performance. In this research, decision-making models are used together with empirical data. One common approach is to derive a decision-making model as the solution of an optimization problem. An objective function $Q_i$ is defined for each option $i$, and the model agent selects the option $i^*$ that maximizes the objective function:
\[ i^* = \argmax_i Q_i. \]
The maximum operation is deterministic and non-differentiable, so for many applications it is replaced by the so-called `softmax' operation, in which option $i$ is chosen with probability
\[ \Pr{i} = \frac{\exp(Q_i)}{\sum_{j=1}^m \exp(Q_j)}. \]
The softmax operation, which we adopt in this paper, is a stochastic, biologically-plausible approximation of the maximum operation~\cite{RSS-AGB:98}. Furthermore, it is differentiable with respect to its argument $Q_i$, which makes it more analytically tractable.

In contexts such as inverse reinforcement learning \cite{SR:98,AYN-SJR:00} and neuroscience \cite{MRN-JIG:13}, a central goal is to understand the decision-making process by finding the objective function values $\{ Q_i \}$ that explain observed decisions. In this paper, we consider this problem in the case that each objective function value $Q_i$ is linear in a set of known variables $\bfx$, i.e., 
\beq \label{eq:linQ}
Q_i = \bftheta^T \bfx_i, \ \ \bftheta, \bfx_i \in \bbR^{n_{obj}}.
\eeq
Models of this form are often used in studies of human decision-making behavior, e.g., \cite{NDD-etal:06, PRM-BKC-JDC:06, JDC-SMM-AJY:07, JG-etal:10}, and are therefore of interest in developing principled methods for human-centered automation. 
By assuming the functional form \eqref{eq:linQ}, we reduce the problem of finding the objective function values to that of learning the vector of parameters $\bftheta$, which we assume to be constant across options and decisions. We call the reduced problem the \emph{parameter estimation problem} for softmax decision-making models with linear objective functions.

The problem of learning the objective function that can explain observed decision-making behavior is relevant for several different disciplines. In the behavioral sciences, it is often of interest to develop models that quantify the various factors that contribute to the decision-making process. Similarly, in engineering, system identification seeks to develop models of dynamic systems that can be used for engineering design. In either case the problem is generally solved in two steps. The first step is to determine which variables affect the process or system in question. In the context of this paper, this is equivalent to determining the variables $\bfx$ in Equation \eqref{eq:linQ}. The second step is to quantify the effect of each variable on the system. This is equivalent to learning the value of the parameters $\bftheta$ in \eqref{eq:linQ}, i.e., solving the parameter estimation problem. We call the two-step process \emph{fitting}. This paper develops an estimator with rigorous performance guarantees for the softmax decision-making model.  This provides a  tool for the second step in the fitting process.

For human-centered automation, a key goal is to develop systems that infer the intuition or the intent of a human operator. One approach is to posit a decision-making model with parameters representing intuition (or intent) and fit the model to observed human choice data. The estimator developed in the present paper makes this possible when applied to an appropriate decision-making model. We demonstrate the estimator using an algorithmic model of human decision-making in a spatial search task, derived in \cite{PR-VS-NEL:14}.  The model, called the stochastic UCL (Upper Credible Limit) model, was derived for multi-armed bandit tasks in a Bayesian setting and was shown to  qualitatively reproduce observed human behavior from experiments.  We use our estimator to infer from these data the human decision-maker's intuition in terms of a set of prior beliefs about the task.  The estimator is applicable to a more general class of decision-making tasks that use a softmax decision-making model.

\begin{figure}
\centering
\includegraphics[width=2.6in]{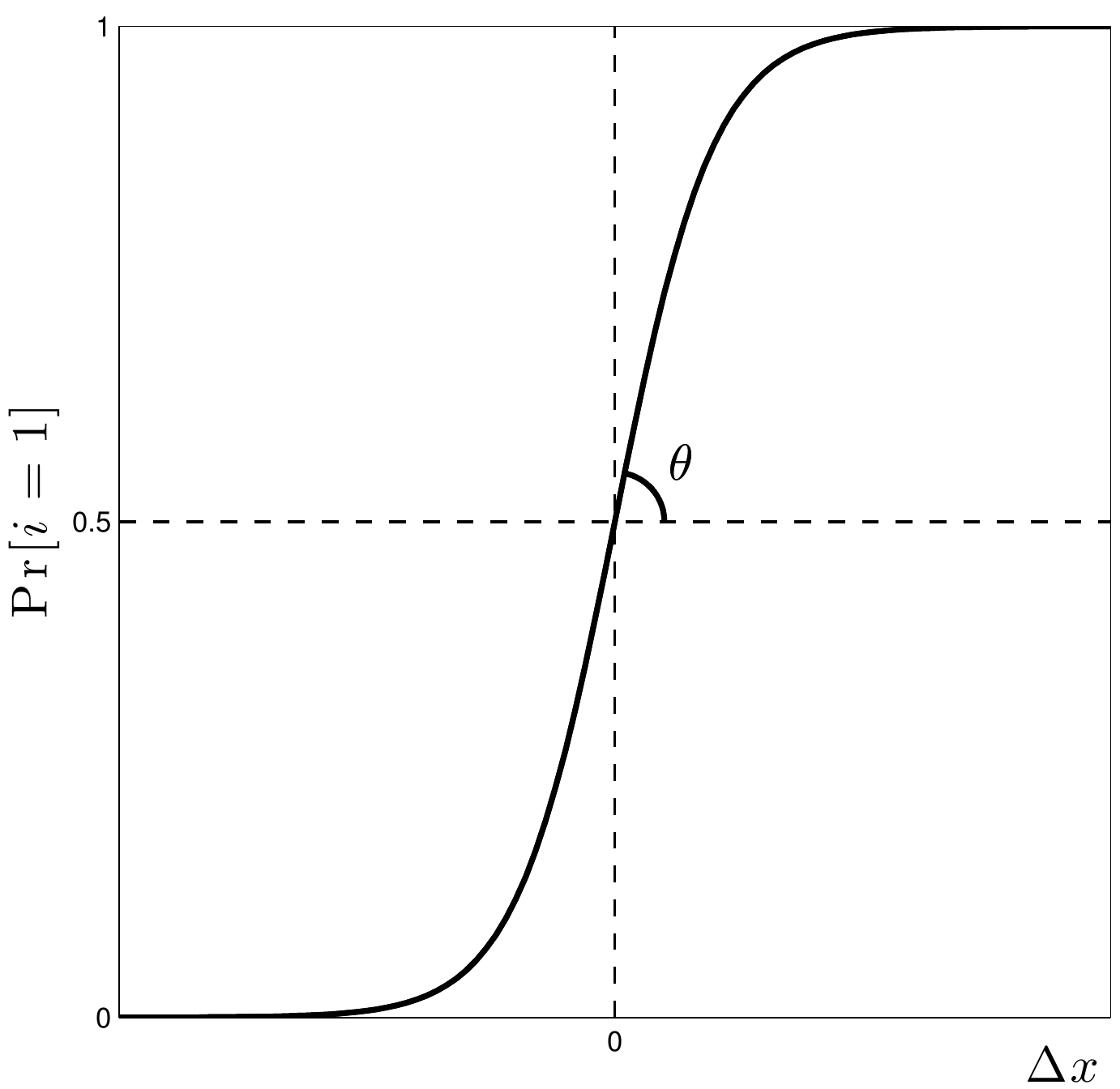}
\caption{The probability \eqref{eq:2Options} from the model \eqref{eq:linQ} with $m=2$ options and a scalar ($n_{obj}=1$) parameter $\theta$. The probability of picking option 1 is a logistic function of $\Delta x = x_1-x_2$ and the sensitivity to $\Delta x$ is controlled by $\theta$, which sets the slope at $\Delta x = 0$.}
\label{fig:2Options}
\end{figure}

As a motivating example of the softmax model, consider the case of deciding between $m=2$ options each with a single ($n_{obj}=1$) known variable  $\bfx_i = x_i$, $i=1,2$, representing the value of the option, and  $\bftheta = \theta$ a scalar. Then the probability of picking option 1 is
\beq \label{eq:2Options}
\Pr{\text{pick option } 1} = \frac{1}{1+ \exp(-\theta (x_1-x_2))}.
\eeq
Figure \ref{fig:2Options} plots the probability \eqref{eq:2Options} as a function of the difference in value of the two options $\Delta x = x_1-x_2$. When the values of the two options are identical, the probability is equal to 0.5 and it increases monotonically with increasing $\Delta x$. The rate of the increase is controlled by $\theta$, which sets the slope of the function at $\Delta x = 0$. Large values of $\theta$ increase the slope and make the choice represented by \eqref{eq:2Options} discriminate between $x_1$ and $x_2$ with more sensitivity, while small values of $\theta$ decrease the slope and make the choice less sensitive to $\Delta x$. Models of this form have been used to study a variety of decision-making tasks \cite{BL-PWG:05,KS-etal:05,NDD-etal:06,AN-etal:12,ARS-etal:12}, where finding the value of $\theta$ that explains a given set of decisions is an important problem.

The parameter estimation problem for softmax decision-making models is related to other problems previously studied in the literature, in particular, the multinomial logistic regression problem \cite{DB:92,BK-LC-MATF-AJH:05} and the conditional log likelihood model learning problem \cite{KG-NAS:10a}. With the linear functional form \eqref{eq:linQ}, the softmax decision-making model and the conditional log likelihood model are formally equivalent, meaning that the parameter estimation problem has been studied in previous work, e.g., \cite{KG-NAS:10a}.  The novelty of the present paper comes in the application of parameter estimators to a formal model of human decision-making (the stochastic UCL model) and its use in quantifying a human subject's intuition about a decision-making task.

The stochastic UCL model for human decision-making in spatial search tasks \cite{PR-VS-NEL:14} is a softmax decision model with an objective function $Q_{UCL}$ that is a nonlinear function of several parameters. We show how $Q_{UCL}$ can be transformed into a linear function of the form \eqref{eq:linQ} by linearizing about a point in parameter space.

We adopt a maximum likelihood approach to parameter estimation. In this framework, the convexity of a model implies asymptotic convergence of estimators and that the estimation problem is a convex optimization problem. The convexity of the conditional log likelihood model is an accepted fact in the natural language processing literature \cite{KG-NAS:10a}, so we do not focus on it here.
We apply the standard optimization algorithms to the stochastic UCL estimation problem and demonstrate our results.

There are two major contributions of this paper. First, we show how to apply standard parameter estimation techniques to the stochastic UCL model, a rigorously-derived model of human choice behavior. Models with a similar softmax functional form are commonly used in the neuroscience literature to model choice behavior and are likely to be widely applicable to the field of human-centered automation. Estimating the parameters of such models provides a method to quantify human intention and intuition in choice tasks. Second, we apply the parameter estimation techniques to empirical human choice data and find statistically significant differences between groups of subjects presented with different experimental conditions.

The remainder of the paper is structured as follows. Section II defines the softmax decision model. Section III defines the parameter estimation problem for the softmax model and reviews convergence results from the literature. Section IV summarizes conditions under which the maximum likelihood estimator converges. Section V provides a numerical example of the estimator. Section VI linearizes  the stochastic UCL model about a nominal parameter $\bar{\bftheta}$ to yield a softmax decision model with a linear objective function, and applies the estimator to simulated data. Section VII applies the linearization procedure to fit the stochastic UCL model to human subject data.  Section VIII concludes.

\section{The Softmax Decision Model}

In this section, we define our notation and the specific softmax decision model for which we derive estimator convergence bounds. We also provide several examples of this model that appear in related literature.

\subsection{Notation}\label{sec:notation}
In the spirit of \cite{BK-LC-MATF-AJH:05}, we set the following notation. We assume we have $n$ observations. For each observation $k$ we have data consisting of $d = m \cdot n_{obj}$ explanatory variables and a response, corresponding to the assignment of one of $m$ classes. Specifically, for each observation $k \in \{1, \ldots, n\}$ we have data $(\bfx^k,\bfy^k)$. For each class $i\in \{1, \ldots, m \}$, we have $n_{obj}$ explanatory variables $\bfx_i^k \in \bbR^{n_{obj}}$. The vector of explanatory variables $\bfx^k \in \bbR^d$ is composed of the concatenation of the $\bfx_i^k$:
\[ \bfx^k = \trans{[ \trans{\bfx_1^k}; \trans{\bfx_2^k}; \cdots \trans{\bfx_m^k} ]}. \]
The response variable $\bfy^k = ( y_1^k, \ldots, y_m^k )^T$ represents the class assignment, where the element $y_i^k = 1$ if the observation corresponds to class $i$ and zero otherwise.

Motivated by models of decision-making \cite{PR-VS-NEL:14}, we consider the following statistical model:
\beq \label{eq:kroneckerModel}
p_i^k(\bftheta) = \Pr{y_i^k = 1 | \bfx^k, \bftheta}  = \frac{\exp\left( \bftheta^{T} \bfx_i^k \right)}{\sum_{j=1}^m \exp \left( \bftheta^{T} \bfx_j^k \right)}
\eeq
for $i \in \{1, \ldots, m \}$, where $\bftheta \in \bbR^{n_{obj}}$ is a weight vector that is the same for all classes. This is the softmax decision-making model with linear objective function \eqref{eq:linQ} introduced above, which has been studied in other literatures under other names. 
In the natural language processing literature, \eqref{eq:kroneckerModel} is known as the \emph{conditional log-likelihood} model, while in the econometrics literature, it is known as the \emph{conditional logit} model \cite{DM:74}.

\subsection{Example softmax decision models}
In this section, we provide several concrete examples of the softmax decision model \eqref{eq:kroneckerModel}. The goal is to make the connection between this functional form and others that appear in the literature.

\begin{example}[Softmax with unknown temperature] \label{ex:softmaxUnkTemp}
A standard decision model in reinforcement learning \cite{RSS-AGB:98} is the so-called softmax action selection rule, which selects an option $i$ with probability
\[ \Pr{i} = \frac{\exp \left( V_i/ \tau \right)}{\sum_{j=1}^n \exp \left( V_j/ \tau \right)}, \]
where $V_i$ is the value associated with option $i$ and $\tau$ is a positive parameter known as the \emph{temperature}. This rule selects options stochastically, preferentially selecting those with higher values. The degree of stochasticity is controlled by the temperature $\tau$. In the limit $\tau \to 0^+$, the rule reduces to the standard maximum and deterministically selects the option with the highest value of $V_i$. In the limit $\tau \to +\infty$, all options are equally probable and the rule selects options according to a uniform distribution.

This model is in the form of \eqref{eq:kroneckerModel} with $n_{obj}=1$. Specifically, assume that the temperature $\tau$ is constant but unknown, and the values $V_i$ are known. Then the two models are identical if we identify
\[ \theta = 1/\tau, \; \;  \bfx_i = V_i. \]
In the reinforcement learning literature, the quantity $1/\tau$ is sometimes known as the \emph{inverse temperature} and referred to by the symbol $\beta$. Our methods allow one to estimate $\theta = 1/\tau = \beta$ from observed choice data.
\end{example}

\begin{example}[Softmax with known cooling schedule form] \label{ex:softmaxLogt}
A slightly more complicated model might let the softmax temperature $\tau$ of Example \ref{ex:softmaxUnkTemp} follow a known functional form, called a \emph{cooling schedule}, that depends on an unknown parameter. For example, in simulated annealing, Mitra~\etal ~\cite{DM-FR-ASV:86} showed that good cooling schedules follow a logarithmic functional form:
\[ \tau(t) = \frac{\nu}{\log t}, \]
where $t$ is the decision time and $\nu > 0$ is a parameter.

If $\nu$ is constant but unknown, this model can be represented in the form of \eqref{eq:kroneckerModel} with $n_{obj}=1$ if we identify
\[ \theta = 1/\nu, \; \;  \bfx_i = V_i \log t.\]
\end{example}

\begin{example}[Softmax $Q$-learning with unknown temperature and learning rate]
According to a simple $Q$-learning model \cite{CW-PD:92}, for each choice time $t$ the agent assigns an expected value $V_i^t$ to each option $i$. The values are initialized to 0 at $t = 1$ and then for each subsequent time, the agent picks option $i_t$, receives reward $r_t$, and updates the value of the chosen option $i_t$ according to 
\[ V_{i_t}^{t+1} = V_{i_t}^t + \alpha \delta_t, \]
where $\alpha \in [0,1]$ is a free parameter called the \emph{learning rate} and $\delta_t = r_t - V_{i_t}^t$ is the \emph{prediction error} at time $t$.

A common model in reinforcement learning \cite{NDD:11} has the agent  make decisions using a softmax rule on the value function $V_i^t$, so the probability of selecting an option $i$ at time $t$ is
\begin{align*}
\Pr{i_t = i} &= \frac{\exp \left( V_i^t/ \tau \right)}{\sum_{j=1}^n \exp \left( V_j^t/ \tau \right)}\\
 &= \frac{\exp \left( V_i^{t-1}/ \tau + \alpha \delta_{t-1} \indicator{i=i_{t-1}} \!/ \tau \right)}{\sum_{j=1}^n \exp \left( V_j^{t-1}/ \tau + \alpha \delta_{t-1} \indicator{i=i_{t-1}} \!/\tau \right)},
 \end{align*}
where $\indicator{\cdot}$ is the indicator function, equal to 1 if its argument is a true statement, and 0 otherwise. Similar models are used in the analysis of fMRI data, e.g. \cite{RCW-YN:13}. If $V_i^{t-1}, V_i^{t-2}$, and $r_t$ are known while $\tau$ and $\alpha$ are unknown, the model is in the form of \eqref{eq:kroneckerModel} with $n_{obj}=2$ if we identify
\[ \bftheta = \begin{bmatrix} 1/\tau; & \alpha/\tau \end{bmatrix}, \;  \; \bfx_i = \begin{bmatrix} V_i^{t-1}; & \delta_{t-1} \indicator{i=i_{t-1}} \end{bmatrix}. \]

If only the initial value $V_i^{t=1} = 0$ is known, then the value function $V_i^t$ becomes a nonlinear function of the parameters $\alpha, \tau$ and the model is not of the form \eqref{eq:kroneckerModel}, although it may be possible to find a transformation that puts it in such a form.
\end{example}

In the following section we define the parameter estimation problem for the softmax model \eqref{eq:kroneckerModel}. We then analyze the problem to develop conditions under which this parameter estimation problem can be solved with provable guarantees about its convergence.

\section{Parameter estimation for softmax decision-making models}
In this section, we define the parameter estimation problem for softmax decision-making models using a likelihood framework, and we review relevant results from the literature. Key to these results is the concept of concavity, which is a property of functions that can guarantee the uniqueness of a maximum. When the likelihood function is concave, the maximum likelihood estimation problem can be solved by off-the-shelf optimization algorithms. Concavity is also central to several results from the econometrics literature that provide conditions under which the estimator is guaranteed to converge asymptotically.

In the optimization literature, it is traditional to consider minimization problems, for which convexity plays the same role as concavity does for maximization problems: a function $f$ is concave if the function $-f$ is convex, and maximizing $f$ is equivalent to minimizing $-f$. Following the literature, we refer to concavity and convexity when discussing results from econometrics and optimization, respectively. We distinguish between two notions of concavity: a function $\map{f}{\bbR^n}{\bbR}$ is \emph{weakly} concave if its Hessian is negative semidefinite, and \emph{strongly} concave if its Hessian is strictly negative definite.

\subsection{The softmax model parameter estimation problem}
In the parameter estimation problem for softmax decision-making models, we wish to estimate the values of $\bftheta$ based on the observed data $(\bfx^k,\bfy^k)$. A standard way to perform parameter estimation is using the maximum likelihood method \cite{SMK:93}. To perform maximum likelihood (ML) estimation of \bftheta, one maximizes the log-likelihood function $\ell(\bftheta)$.
\begin{prob} \label{prob:mlProblem}
The maximum likelihood \emph{parameter estimation problem} for the softmax decision model \eqref{eq:kroneckerModel} is the optimization problem
\beq \label{eq:mlProblem}
\hat{\bftheta}_{ML} = \argmax_{\bftheta} \ell(\bftheta),
\eeq
where $\ell(\bftheta)$ is the logarithm of the likelihood function of the model \eqref{eq:kroneckerModel}, defined as
\begin{align} \label{eq:logLikelihood}
\ell&(\bftheta) = \sum_{k=1}^n \log \Pr{\bfy^k | \bfx^k, \bftheta}\\
&= \sum_{k=1}^n \left[ \sum_{i=1}^m y_i^k \bftheta^{T} \bfx_i^k - \log \sum_{i=1}^m \exp \left( \bftheta^{T} \bfx_i^k \right) \right].\nonumber
\end{align}
The ML estimate $\hat{\bftheta}_{ML}$ can be interpreted as the parameter value $\bftheta$ that makes the observed data most likely under the given model.
\end{prob}

A prior on \bftheta\ can be incorporated by adopting a maximum a posteriori (MAP) estimate,
\beq
\hat{\bftheta}_{MAP} = \argmax_{\bftheta} L(\bftheta) = \argmax_{\bftheta} [\ell(\bftheta) + \log p(\bftheta) ], \label{eq:mapProblem}
\eeq
with $p(\bftheta)$ being the prior on \bftheta. The MAP estimate penalizes ML estimates that are considered unlikely under the prior.

\subsection{Asymptotic behavior of the ML estimator} \label{sec:mlAsymptotics}
The ML estimator $\hat{\bftheta}_{ML}$ solves the estimation problem in the frequentist framework, which posits that there is a true value $\bftheta_0$ of the parameters that we attempt to recover from analyzing the given data. In this framework, natural questions to be asked are 1) does $\hat{\bftheta}_{ML} \to \bftheta_0$ as the number of observations $n$ grows, and 2) how dispersed is the difference $\hat{\bftheta}_{ML} - \bftheta_0$? These questions have been studied in the econometrics literature, for which \cite{WKN-DLM:94} is a standard reference. The remainder of this section summarizes the relevant results from \cite{WKN-DLM:94}. The answers to these two questions depend on two properties of the model, identification and concavity, defined as follows.

\begin{defn}[Identification]
A statistical model with likelihood function $\map{\ell}{\bbR^q}{\bbR}$ and observed data $\bfx$ is said to be \emph{identified} if, for all $\bftheta, \bftheta_0 \in \bbR^q$, 
\[ \bftheta \neq \bftheta_0 \Rightarrow \ell(\bftheta_0 ; \bfx) \neq \ell(\bftheta ; \bfx). \]
\end{defn}

\begin{defn}[Concavity]
A statistical model with likelihood function $\map{\ell}{\bbR^q}{\bbR}$ is said to be \emph{concave} if $\ell(\bftheta ; \bfx)$ is strictly concave in $\bftheta$.
\end{defn}

If a model with likelihood function $\ell$ is identified and concave (see \cite[Theorem 2.7]{WKN-DLM:94} for details), the answer to question 1) is yes. These two properties imply that the true value $\bftheta_0$ of the parameter is the unique maximum of the expected value of the log-likelihood $\ell(\bftheta)$.

Concavity and identification can depend on both the functional form of $\ell(\bftheta)$ and the observed data $\bfx$. As an example of how the identification property may fail due to data, consider the model \eqref{eq:kroneckerModel} with $\bfx_i$ being the zero vector for each $i$. In this case, $\Pr{y_i = 1 | \bfx, \bftheta} = 1/m$ for each $i$ independent of $\bftheta$ and the estimation procedure will be unable to distinguish among the possible parameter values. In the following section, we derive conditions on the data that ensure identification. These conditions also ensure that $\ell(\bftheta)$ is strictly concave and provide guidelines for the design of experiments for estimating $\bftheta$.

The answer to question 2) is that, under mild regularity conditions, the distribution of $\hat{\bftheta}_{ML}$ approaches a normal distribution as the number of samples $n$ grows. In particular, the following limit holds:
\beq \label{eq:asymptoticNormality}
\hat{\bftheta}_{ML} \stackrel{d}{\to} \mcN(\bftheta_0,\bfJ^{-1}/n),
\eeq
where $\stackrel{d}{\to}$ signifies a limit in distribution as $n \to \infty$ and $\bfJ = -\E{\bfH(\bftheta_0)}$ is the negative of the expected value of the Hessian of $\ell(\bftheta)$ with respect to $\bftheta$. See \cite[Chapter 9]{ASG:91} for more details about the concept of a limit in distribution and see \cite[Theorem 3.3]{WKN-DLM:94} for full details of the conditions under which \eqref{eq:asymptoticNormality} holds. In practice one uses $\hat{\bfJ} = - \bfH(\hat{\bftheta}_{ML})/n$ as an estimate of $\bfJ$. This permits construction of standard frequentist analysis tools, such as confidence intervals for the parameter estimates and hypothesis tests. The estimate $\hat{\bftheta}_{ML}$ is efficient in the sense that it obeys the Cram\'er-Rao lower bound \cite{SMK:93} on the variance of estimators $\hat{\bftheta}$, so no other unbiased estimator can have lower variance than $\hat{\bftheta}_{ML}$.

\section{Analysis of the maximum likelihood estimator for softmax decision models with linear objective functions}
In this section we present conditions under which the model \eqref{eq:kroneckerModel} is identified and concave. These conditions imply that the ML estimator $\hat{\bftheta}_{ML}$ converges and that the ML optimization problem \eqref{eq:mlProblem} is convex. The concavity of the model is an accepted fact in the natural language processing literature \cite{KG-NAS:10a};  we summarize the result in Theorem \ref{thm:mlConvergence}.

\subsection{Asymptotic and finite-sample behavior}
Recall from Section \ref{sec:mlAsymptotics} that two properties that guarantee asymptotic convergence of the ML estimator are identification and concavity. Whether or not the model \eqref{eq:kroneckerModel} satisfies these properties can be a function of the data $\bfx^k, k \in \{1, \ldots, n \}$. Recall our example where $\bfx_i^k = \mathbf{0}$ for each $i$ and $k$. In this case the probability $\Pr{y_i^k | \bfx^k, \bftheta} = 1/m$ for each $i$ and $k$ for all values of $\bftheta$ and the likelihood function $\ell(\bftheta)$ is flat, so neither identification nor concavity is satisfied.

However, a sufficient condition for identification is as follows. Define the $n_{obj} \times m$ matrix $\mathbf{X}^k$ by transforming the explanatory variable $\bfx^k$ of a single observation $k$:
\beq \label{eq:bfX}
\mathbf{X^k} = \begin{bmatrix} \bfx^k_1 & \bfx^k_2 & \cdots & \bfx^k_{m-1} & \mathbf{0} \end{bmatrix}.
\eeq
Note that $\bfX^k \trans{\bfX^k} = \sumOp(\trans{\bfx^k} \!\! \bfx^k)$. Considering $\bfX^k$ as a random variable, the following lemma ensures identification.
\begin{lem} \label{lem:id}
Let $\bfx$ be the explanatory variable for an arbitrary observation and let $\bfX$ be the transformation of $\bfx$ defined in \eqref{eq:bfX}. If the second-moment matrix $\E{\bfX \trans{\bfX}}$ exists and is positive definite, then the model \eqref{eq:kroneckerModel} is identified.
\end{lem}

\begin{IEEEproof}
The probability of choosing an option $i$ under the model \eqref{eq:kroneckerModel} is a monotonic function of the objective value $Q_i$, so it suffices to show that the data provides a one-to-one mapping between the parameter vector $\bftheta$ and the objective values $Q_1, \ldots, Q_m$.

Let $\bftheta, \bftheta^{\prime} \in \bbR^{n_{obj}}$ and define the vectors of objective function values $\mathbf{Q} = \bftheta^T \bfX$ and $\mathbf{Q}^{\prime} = \bftheta^{\prime T} \bfX$. Define $\Delta \mathbf{Q} = \mathbf{Q}-\mathbf{Q}^{\prime} = (\bftheta-\bftheta^{\prime})^T \bfX \in \bbR^m$. The magnitude of $\Delta \mathbf{Q}$ satisfies $\E{\| \Delta \mathbf{Q} \|^2} = (\bftheta-\bftheta^{\prime})^T \E{\bfX \bfX^T} (\bftheta-\bftheta^{\prime})$. Then by the assumption that $\E{\bfX \bfX^T}$ is positive definite, $\E{\| \Delta \mathbf{Q} \|^2} = 0$ implies $(\bftheta - \bftheta^{\prime}) = 0$, so $\bftheta = \bftheta^{\prime}$ and $\mathbf{Q} = \mathbf{Q}^{\prime}.$ Therefore the mapping between the parameters $\bftheta$ and the objective values $Q_1, \ldots, Q_m$ is one-to-one, which implies that $\ell(\bftheta|\bfx,\bfy) \neq \ell(\bftheta^{\prime}|\bfx,\bfy)$ for $\bftheta \neq \bftheta^{\prime}$ and the model is identified.
\end{IEEEproof}

The condition of Lemma \ref{lem:id} is given in terms of an expectation, but in practice one has a given sample of data. In this case the expectation can be replaced by the sample average. Specifically, define $\bfX^k$ for each observation $k$ following \eqref{eq:bfX}. Then $\E{\bfX \bfX^T}$ is estimated by
\[ \E{\bfX \bfX^T} \approx \frac{1}{n} \sum_{j=1}^n \bfX^k \bfX^{kT}. \]
If this sample average is positive definite, then the model is identified. For the sample average to be positive definite it must be full rank $=n_{obj}$, and each observation $k$ can add at most $m$ to the rank. Therefore, the following inequality must be satisfied for the model to be identified:
\[ m \cdot n \geq n_{obj}. \]
This gives a lower bound $n \geq \lceil n_{obj}/m \rceil$ on the minimum number of observations required for identification. For most applications, the number of options $m$ will be larger than the number of parameters $n_{obj}$, so the lower bound is trivial,.  However, for cases with large number of parameters the bound can be useful for experimental design.

The following theorem summarizes the conditions under which the ML estimator \eqref{eq:mlProblem} converges.
\begin{thm}[Convergence of the ML estimator] \label{thm:mlConvergence}
Let $\bfX^k$ be defined as in \eqref{eq:bfX}. If the second-moment matrix 
\[ \frac{1}{n} \sum_{k=1}^n \bfX^k \trans{\bfX^{k}} \]
exists and is positive definite, then 
\begin{enumerate}
\item The ML optimization problem \eqref{eq:mlProblem} is convex.
\item The ML estimator $\hat{\bftheta}_{ML}$ for the model \eqref{eq:kroneckerModel} is asymptotically approximately distributed as
\beq \label{eq:thetaMLAsymptotic}
\hat{\bftheta}_{ML} \sim \mcN(\bftheta_0,\hat{\bfJ}^{-1}/n),
\eeq
where $\hat{\bfJ} = - \bfH(\hat{\bftheta}_{ML})/n$ is the empirical mean Hessian of the likelihood function evaluated at the estimated parameter value.
\end{enumerate}
\end{thm}
\begin{IEEEproof}
See \cite{ReverdyThesis} and \cite{DM:74}.
\end{IEEEproof}

Theorem \ref{thm:mlConvergence} proves convergence of the parameter estimate $\hat{\bftheta}_{ML}$ and provides its asymptotic distribution. This distribution can be used to formulate frequentist confidence intervals for the parameter estimate $\hat{\bftheta}_{ML}$. Furthermore, the theorem proves that the optimization problem \eqref{eq:mlProblem} is convex, which allows us to solve it using off-the-shelf optimization algorithms. In the following, we use the phrase \emph{the estimator} to refer to the procedure of using an off-the-shelf convex optimization algorithm to solve the maximum likelihood problem \eqref{eq:mlProblem}. We use the phrase \emph{the estimate} to refer to the solution $\hat{\bftheta}$ of \eqref{eq:mlProblem} thus obtained.

\section{Numerical examples}
In this section we present several numerical examples to demonstrate the theory developed in the previous sections for solving the parameter estimation problem \eqref{eq:mlProblem}.

\subsection{Scalar parameter}
First, we consider model \eqref{eq:kroneckerModel} with $m=10$ options and a scalar parameter $\bftheta = \theta = \theta_0$ that we wish to estimate. This could correspond to a  decision-maker  choosing among ten options using a softmax model with unknown constant inverse temperature $\theta = \theta_0$, as in Example \ref{ex:softmaxUnkTemp}. Alternatively, it could correspond to a temperature varying with  observation number $k= 1, \ldots, n$ according to a known function with a single unknown parameter $\theta = \theta_0$, e.g., $\tau_k = \theta/\log k$, as in Example \ref{ex:softmaxLogt}. In this case the $\log k$ term can be absorbed into the explanatory variables and we proceed as in the first case.

\begin{figure}
\centering
\includegraphics[width=3in]{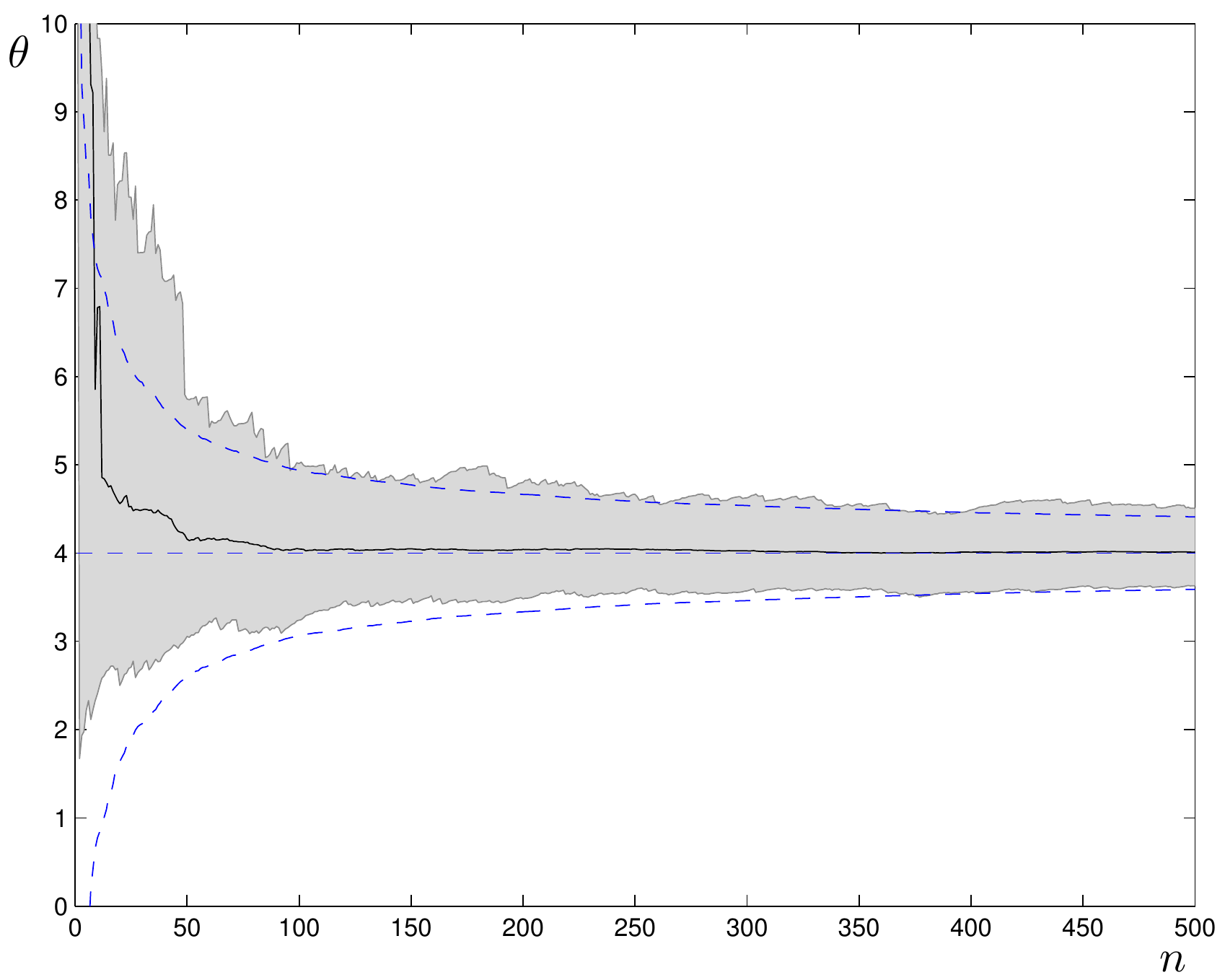}
\caption{Scalar parameter estimation example. Illustration of the convergence of parameter estimates to the asymptotic normal distribution \eqref{eq:thetaMLAsymptotic} as the number of observations $n$ grows. The dashed lines show the true value of the scalar parameter $\theta_0=4$ and the accompanying 95\% confidence intervals implied by the asymptotic normal distribution \eqref{eq:thetaMLAsymptotic}. For each value of $n$, an ensemble of 100 parameter estimates was formed by repeatedly simulating the data $\bfy$ while holding the explanatory variables $\bfx$ fixed, and using the estimator to compute the value of the parameter. The solid black line shows the mean parameter estimate and the shaded region the empirical 95\% confidence interval.}
\label{fig:scalarConvergence}
\end{figure}

Figure~\ref{fig:scalarConvergence} shows results of applying the estimator to simulated data. For every $k$, when an observation was taken and a decision made, the model was simulated 100 times.  For each of the 100 simulations,  the estimator was applied to estimate the parameter $\theta$ based on the first $k$ observations.  Running 100 simulations made it possible to examine convergence of the estimate in distribution.  Figure~\ref{fig:scalarConvergence} illustrates how the estimates converge in distribution to the normal distribution \eqref{eq:thetaMLAsymptotic} as the number of observations $n$ increases.  For the simulations, the explanatory variables $\bfx^k$ were drawn from a standard Gaussian distribution $\mcN(0,1)$ (mean zero and unit variance), and the response variables $\bfy^k$ were drawn according to probability distribution \eqref{eq:kroneckerModel} conditional on $\bfx^k$ and $\theta_0 = 4$. The estimates were computed by solving the optimization problem \eqref{eq:mlProblem} using a BFGS quasi-Newton algorithm \cite{CGB:70,RF:70,DG:70,DFS:70} (Matlab function \verb.fminunc.). Theorem \ref{thm:mlConvergence} guarantees that the optimization problem is convex, so the algorithm will converge.

The convergence behavior can be seen in Figure \ref{fig:scalarConvergence} by observing the mean parameter estimate $\hat\theta_{ML}$ as well as its confidence intervals. The mean parameter estimate $\hat\theta_{ML}$, represented by the solid black line, converges to the true parameter value $\theta_0 = 4$, represented by the horizontal dashed line. However, Theorem \ref{thm:mlConvergence} guarantees convergence in distribution, which is a stronger result. To illustrate this behavior we plot 95\% confidence intervals for both the empirical distribution of estimates $\hat{\theta}_{ML}$ and the asymptotic distribution \eqref{eq:thetaMLAsymptotic}, computed from the ensemble of 100 parameter estimates. For values of $n$ greater than 100, the two intervals overlap closely, showing that the distribution of estimates has converged. Importantly, this shows that statistical hypothesis tests based on the asymptotic distribution \eqref{eq:thetaMLAsymptotic} will be accurate.

For small amounts of data, i.e., $n < 50$, the mean parameter estimate is biased above the true value. The bias is due to an insufficient amount of data being used in the estimation procedure, and the direction of the bias can be explained as follows. Larger values of the parameter $\theta$ correspond to more deterministic choice behavior. When $\theta_0 > 0$, for any given choice, the model is more likely to pick the option with a larger objective value, resulting in a parameter estimate that is biased upwards. This bias can be seen in Figure~\ref{fig:vectorConvergence} as well, which treats a case with a vector parameter.

\subsection{Vector parameter}
\begin{figure}
\centering
\includegraphics[width=3in]{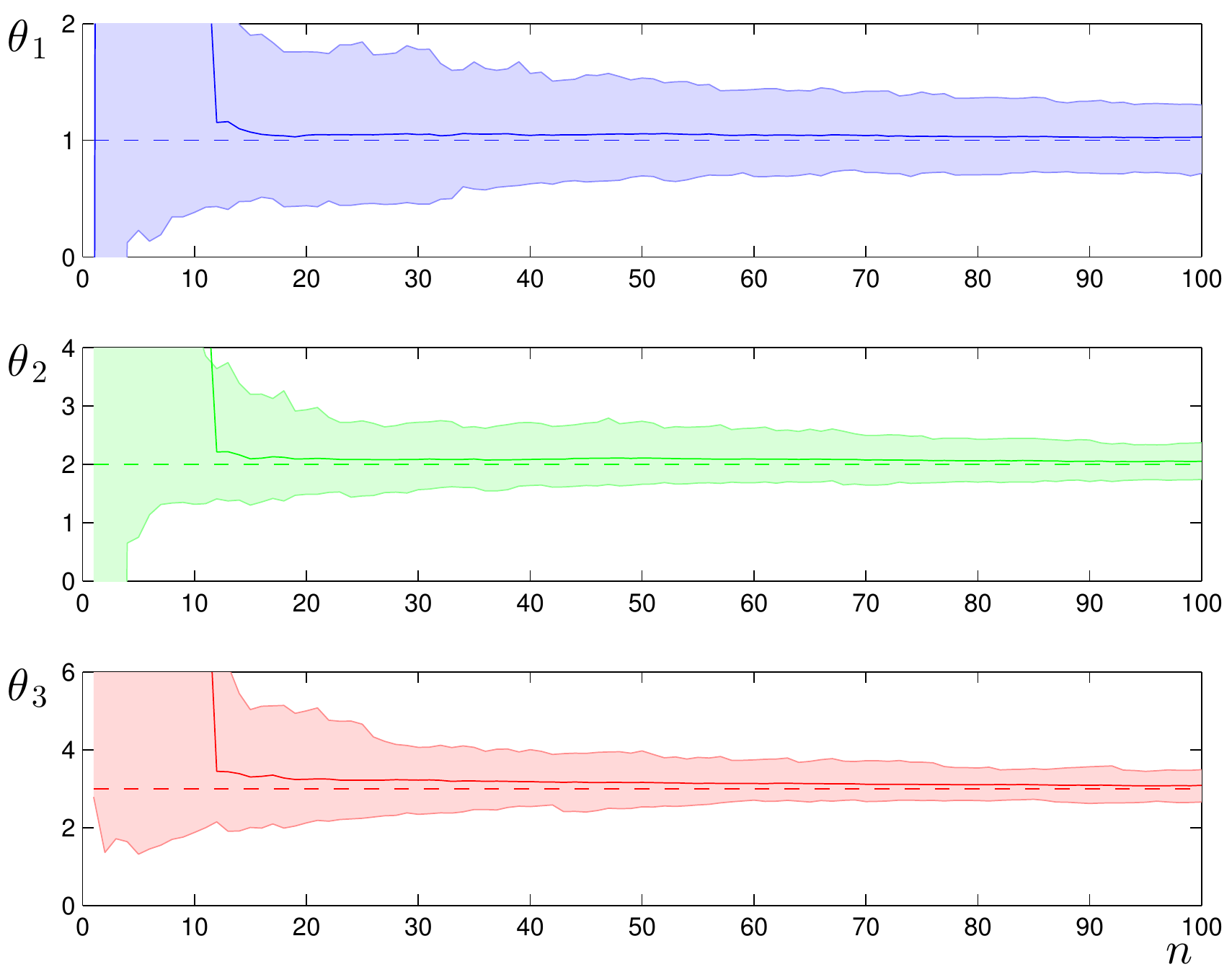}
\caption{Vector parameter estimation example. Illustration of the convergence of parameter estimate to the asymptotic normal distribution \eqref{eq:thetaMLAsymptotic} as the number of observations $n$ grows. The dashed  lines show the true value of each element of the vector parameter $\bftheta_0=[1,2,3]^T$. For each value of $n$, an ensemble of 100 parameter estimates was formed by repeatedly simulating the data $\bfy$ while holding the explanatory variables $\bfx$ fixed, and using the estimator to compute the value of the parameter. The solid lines show the mean parameter estimate and the shaded regions the empirical 95\% confidence interval.}
\label{fig:vectorConvergence}
\end{figure}

Next, we consider the model \eqref{eq:kroneckerModel} with $m = 100$ options and a vector parameter $\bftheta = \bftheta_0$ with $n_{obj}=3$ elements that we wish to estimate.  Figure~\ref{fig:vectorConvergence} shows results of applying the estimator to simulated data in this vector parameter estimation example. As in the scalar parameter estimation case above, the model was simulated 100 times for every $k = 1, \ldots, n$.   Figure~\ref{fig:vectorConvergence} shows how the estimate converges to the true value $\bftheta_0 $ as the total number of observations $n$ increases.  The explanatory variables $\bfx^k$ were drawn according to independent standard Gaussian distributions, and the response variables $\bfy^k$ drawn according to the model \eqref{eq:kroneckerModel} conditional on $\bfx^k$ and true vector parameter value $\bftheta_0= [1,2,3]^T$.  The estimates were computed as in the scalar case.

The convergence behavior can be seen in Figure \ref{fig:vectorConvergence} by observing the mean parameter estimate $\hat\bftheta_{ML}$ as well as its confidence intervals.
For each of the three parameters $\theta_i$, $i=1, 2, 3$, the corresponding mean parameter estimate $\left(\hat\theta_{ML}\right)_i$, represented as a solid line, converges to the true parameter value $(\bftheta_0)_i$, represented by a horizontal dashed line. The shaded regions represent the empirical 95\% confidence interval around the corresponding mean value, computed from the ensemble of 100 parameter estimates. For clarity, we omit the confidence intervals implied by the asymptotic normal distribution \eqref{eq:thetaMLAsymptotic} from the figure, but the behavior is similar to that shown in Figure~\ref{fig:scalarConvergence}.

There is an upwards bias in the parameter estimates for small numbers of observations $n$, as in Figure \ref{fig:scalarConvergence}. The width of the confidence intervals for the three parameters scales roughly with their true value $(\bftheta_0)_i$. This behavior can be seen in the figures in the next section as well.

\section{Application to nonlinear objective functions using linearization} \label{sec:linearization}
The development up to this point for addressing the parameter estimation problem \eqref{eq:mlProblem} has assumed that the objective function takes the linear form \eqref{eq:linQ}. However, many relevant objective functions are nonlinear functions of the unknown parameter $\bftheta$.  One approach is to linearize the nonlinear objective function about a nominal parameter value, and then apply the estimator to the linearized objective function.  We apply this approach to the nonlinear objective function from the stochastic UCL algorithm \cite{PR-VS-NEL:14}, an algorithm for human decision-making in  multi-armed bandit tasks in a Bayesian setting, and show how its parameters can be estimated.

\subsection{The multi-armed bandit problem}
The multi-armed bandit problem, introduced by Robbins \cite{HR:52} is a sequential decision problem which consists of a set of $N$ options (each option is also called an \emph{arm} in analogy with the lever of a slot machine). Each option $i \in \{1, \ldots, N\}$, has an associated probability distribution $p_i$ with mean $m_i$, unknown to the agent solving the problem. At each sequential decision time $t \in \{1, \dots, T\}$ the agent picks an arm $i_t$ and receives a stochastic reward $r_t \sim p_{i_t}(r)$ drawn from the probability distribution associated with that arm. This is a special case of the notation introduced in Section \ref{sec:notation}, with $m = N$ options indexed by $i$ and $n = T$ decisions indexed by $t$. The agent's objective is to maximize the expected value of the cumulative rewards received from the $T$ decisions:
\[ \max_{\{i_t\}} J, \;\;\; \ J = \E{\sum_{t=1}^T r_t} = \sum_{t=1}^T m_{i_t}. \]

Each choice of $i_t$ is made conditional on the information available to the agent at time $t$. If the mean rewards $m_i$ were known to the agent, the optimal policy would be trivial: pick arm $i_t \in \argmax_i m_i$ for each $t$. However, since the mean rewards are unknown, the agent must simultaneously select arms where the reward value is uncertain to gain information about the rewards and preferentially select arms with high rewards to accumulate reward. The tension between selecting arms with uncertain (but possibly high) rewards  and selecting arms that appear to have high rewards based on current information is known as the \emph{explore-exploit} tradeoff.  This tradeoff is common to a variety of problems in machine learning and adaptive control.

The multi-armed bandit problem is the subject of active research in machine learning as well as in neuroscience. In \cite{PR-VS-NEL:14}, we showed that a significant fraction of human subjects exhibited excellent performance in solving a multi-armed bandit problem, even outperforming algorithms known to have optimal performance in some cases. We attributed this good performance to the human subjects' having good priors on the structure of the rewards $m_i$, and we designed the stochastic UCL algorithm as a model of human performance to capture this dependence on priors. Estimating the parameters of this model from observations of a human solving the multi-armed bandit task would allow a machine to learn the human's belief priors.  This could in turn facilitate a design of a human-machine system that could achieve better performance than either the human or the machine could on its own.

\subsection{The stochastic UCL algorithm}
The stochastic UCL algorithm \cite{PR-VS-NEL:14} is designed to solve multi-armed bandit problems with Gaussian rewards, i.e., where the reward distribution $p_i(r) = \mcN(m_i,\sigma_s^2)$ is Gaussian with unknown mean $m_i$ and known variance $\sigma_s^2$. The algorithm consists of two parts: Bayesian inference that maintains the agent's belief state and a softmax decision model that uses an objective function $Q$ that depends on the belief state. Both the inference and the decision parts introduce nonlinear dependencies on the parameters of the algorithm.

As a model of human behavior, the stochastic UCL algorithm assumes that the agent's prior distribution of $\mathbf{m}$ (i.e., the agent's initial beliefs about the mean reward values $\bfm$ and their covariance) is multivariate Gaussian with mean $\boldsymbol{\mu_0}$ and covariance $\Sigma_0$:
\[ \bfm \sim \mcN(\boldsymbol{\mu_0},\Sigma_0),\]
where $\boldsymbol{\mu_0} \in \mathbb{R}^{N}$ and $\Sigma_0 \in \mathbb{R}^{N \times N}$ is a positive-definite matrix.

In \cite{PR-VS-NEL:14} we use a minimal set of three parameters to specify $(\bfmu_0,\Sigma_0)$. For the mean we use a uniform prior $\bfmu_0 = \mu_0 \mathbf{1}_N,$ where $\mu_0 \in \bbR$ is a single parameter that encodes the agent's belief about the mean value of the rewards and $\oneVec{N}$ is the vector with each element equal to 1. For the problems considered in \cite{PR-VS-NEL:14}, the arms are spatially embedded with each arm at a different location in space (see Figure \ref{fig:surfaceProfiles} in the next section). It is reasonable to assume that arms that are spatially close will have similar mean rewards. Therefore, for the covariance $\Sigma_0$ we set $\Sigma_0 = \sigma_0^2 \Sigma$ where $\Sigma$ represents a prior that is exponential in distance, i.e., each element has the form
\beq \label{eq:exponentialPrior}
\Sigma_{ij} = \exp(-\| z_i - z_j \|/\lambda),
\eeq
where $z_i$ is the location of arm $i$ and $\lambda \geq 0$ is the correlation length scale. The parameter $\sigma_0 \geq 0$ can be interpreted as a confidence parameter, with $\sigma_0 = 0$ representing absolute confidence in the beliefs about the mean $\bfmu_0$, and $\sigma_0 = +\infty$ representing complete lack of confidence.

With this prior, the posterior distribution is also Gaussian, so the Bayesian optimal inference algorithm is linear and can be written down as follows. At each time $t$, the agent selects option $i_t$ and receives a reward $r_t$. Let $\bfr^t$ be the $t \times 1$ vector composed of the $r_t$. Let $n_{i}^t$ be the number of times the agent has selected option $i$ up to time $t$, let $\bar{m}_{i}^t$ be the empirical mean reward observed for option $i$, and let $\bfn^t$ and $\bar{\bfm}^t$ be the vectors composed of the $n_i^t$ and $\bar{m}_i^t$, respectively. For each time $t$, define the precision matrix $\Lambda_t = \Sigma_t^{-1}$. Then the belief state at time $t$ is \cite[Theorem 10.3]{SMK:93}
\begin{align}
\Lambda_t &= \frac{\diag(\bfn^t)}{\sigma_s^2} + \Lambda_0, \ \Sigma_t = \Lambda_t^{-1} \label{eq:GenInferenceSigma}\\
\bfmu_t &= \bfmu_0 + \Sigma_0 H_t^T \! \left( H_t \Sigma_0 H_t^T + \sigma_s^2 I_t \right)^{-1} \! (\bfr^t - H_t \bfmu_0), \label{eq:GenInferenceMu}
\end{align}
where $H_t$ is the $t \times N$ observation matrix with $H_t(t,j) = 1$ if $i_t = j$ and zero otherwise, and $I_t$ is the $t$-dimensional identity matrix.

Based on the belief state $(\bfmu_t,\Sigma_t)$, the stochastic UCL algorithm chooses arm $i_t$ with probability
\beq \label{eq:prUCL}
\Pr{i_t = i | \tilde{Q}^t, \upsilon_t} = \frac{\exp(\tilde{Q}_i^t/\upsilon_t)}{\sum_{j=1}^N \exp(\tilde{Q}_i^t/\upsilon_t)},
\eeq
where $\tilde{Q}_i^t$ is the heuristic function value for arm $i$ at time $t$ and $\upsilon_t$ is the temperature corresponding to the cooling schedule at time $t$. The cooling schedule is assumed to take the form $\upsilon_t = \nu/\log t$, $\nu$ a constant, so the probabilities \eqref{eq:prUCL} become
\beq \label{eq:softmaxUCLPr}
\Pr{i_t = i | \tilde{Q}^t, \nu} = \frac{\exp((\tilde{Q}_i^t \log t)/\nu)}{\sum_{j=1}^N \exp((\tilde{Q}_i^t \log t)/\nu)}.
\eeq
The heuristic function is
\beq \label{eq:Qit}
\tilde{Q}_i^t = \mu_i^t + \sigma_i^t \Phi^{-1}(1-\alpha_t),
\eeq
where $\mu_i^t = \left(\bfmu_t\right)_i$ is the posterior mean reward of arm $i$ at time $t$ and $\sigma_i^t = \sqrt{\left(\Sigma_t\right)_{ii}}$ its associated standard deviation. The quantity $\Phi^{-1}(\cdot)$ is the inverse cumulative distribution function of the standard normal distribution and $\alpha_t = 1/\sqrt{2 \pi e}t$ is a decreasing function of time.

This is a softmax decision model with unknown parameters $(\mu_0,\sigma_0,\lambda,\nu)$, but it is not yet in the form \eqref{eq:kroneckerModel} since the quantity $(\tilde{Q}_i^t \log t)/\nu$ is a nonlinear function of the parameters. However, we can locally approximate $(\tilde{Q}_i^t \log t)/\nu$ with a linear function by linearizing about a nominal value of the prior. By estimating the parameter values of the linearized model, we can recover the parameters of the original nonlinear model \eqref{eq:softmaxUCLPr} near the nominal prior.

\subsection{Linearization}
Let $\delta_0^2 = \sigma_s^2/\sigma_0^2$ be the relative precision of a reward measurement compared to the certainty of the prior. Fix a nominal prior with parameter values $(\muBar, \deltaBar,\lambda)$ and consider small deviations $\deltaMu$ and $\deltaDelta$ about $\muBar$ and $\deltaBar$, respectively:
\[ \mu_0 = \muBar + \deltaMu, \ \delta_0^2 = \deltaBar + \deltaDelta. \]
In the case that the true value of $\lambda$ is unknown, this method is easily generalized to include deviations in $\lambda$, but for simplicity of exposition we consider it fixed. Recall that the covariance prior is $\Sigma_0 = \sigma_0^2 \Sigma$, where $\Sigma$ is defined by  \eqref{eq:exponentialPrior}, and  its inverse is denoted by $\Lambda = \Sigma^{-1}$.

In terms of the nominal value $\bar{\delta}_0^2$,  \eqref{eq:GenInferenceSigma} becomes
\[ \Lambda_t = \frac{1}{\sigma_s^2} \left( \diag(\bfn^t) + \deltaBar \Lambda + \deltaDelta \Lambda \right). \]
Therefore, to first order in \deltaDelta, $\Sigma_t$ is given by
\beq \label{eq:SigmatLinear}
\Sigma_t = \sigma_s^2 A_t^{-1} - \sigma_s^2 A_t^{-1}BA_t^{-1} \deltaDelta + \bigO{\deltaDelta^2},
\eeq
where $A_t = \deltaBar \Lambda + \diag(\bfn^t)$ and $B = \Lambda = \Sigma^{-1}$. Expanding the square root in the following, we get 
\beq
\sigma_i^t = \sqrt{\left(\Sigma_t\right)_{ii}} = \sqrt{c_i^t} - \frac{d_i^t}{2 \sqrt{c_i^t}} \deltaDelta + \bigO{\deltaDelta^2},
\eeq
where $c_i^t$ is the $i^{th}$ element on the diagonal of $C_t = \sigma_s^2 A_t^{-1}$ and $d_i^t$ is the $i^{th}$ element on the diagonal of $D_t = \sigma_s^2 A_t^{-1}BA_t^{-1}$. The standard deviation $\sigma_i^t$ must be nonnegative, which implies an upper bound on $\deltaDelta$. Similarly, $\delta_0^2$ must be nonnegative, which implies a lower bound on $\deltaDelta$, which is already assumed to be small. The implied bounds on $\deltaDelta$ are
\[ - \deltaBar = - \frac{\sigma_s^2}{\bar{\sigma}_0^2} \leq \deltaDelta \leq \frac{2 c_{i}^t}{d_{i}^t}, \]
which, together with the requirement that $\deltaDelta$ be small with respect to $\deltaBar$, gives a bound on the values of $\deltaDelta$ for which the linearization is valid.

Similarly, the expression \eqref{eq:GenInferenceMu} for $\bfmu_t$ becomes
\beq \label{eq:mutLinear}
\bfmu_t = E_t + F_t \deltaMu + G_t \deltaDelta + \bigO{\Delta^2},
\eeq
where $\Delta^2$ denotes second-order terms in the deviation variables $\deltaDelta$ and $\deltaMu$, and $E_t, F_t$, and $G_t$ are the $N \times 1$ vectors
\begin{align}
E_t &= \muBar \oneVec{N} + \frac{\Sigma H_t^T}{\deltaBar} \! \! \left(I_t - H_t A_t^{-1} H_t^T \right) \! (\bar{\bfm}^t - H_t \muBar \oneVec{N})\\
F_t &= \oneVec{N} - \frac{\Sigma H_t^T}{\deltaBar} \! \! \left(I_t - H_t A_t^{-1} H_t^T \right) \! H_t \oneVec{N}\\
G_t &= -A_t^{-1} B A_t^{-1}(H_t^T \bfm^t - \bfn^t \muBar).
\end{align}

Define $e_i^t, f_i^t$, and $g_i^t$ as the $i^{th}$ components of $E_t, F_t,$ and $G_t$, respectively. Then the linearized heuristic is
\beq \label{eq:linearizedHeuristic}
\frac{\tilde{Q}_i^t \log t}{\nu} \approx Q_i^t = \bftheta^T \bfx_i^t = \theta_1 x_{i,1}^t + \theta_2 x_{i,2}^t + \theta_3 x_{i,3}^t ,
\eeq
where the parameters $\bftheta$ are defined by
\beq \label{eq:linParamDef}
\theta_1 = \frac{1}{\nu},  \;\;\;  \theta_2 = \frac{\deltaMu}{\nu}, \;\;\;  \theta_3 = \frac{\deltaDelta}{\nu}
\eeq
and the explanatory variables $\bfx_i^t$ are defined as
\begin{align}
x_{i,1}^t &= \left(e_i^t + \sqrt{c_i^t} \Phi^{-1}(1-\alpha_t) \right)\log t\\
x_{i,2}^t &= f_i^t \log t\\
x_{i,3}^t &= \left(g_i^t - \frac{d_i^t}{2\sqrt{c_i^t}} \Phi^{-1}(1-\alpha_t) \right) \log t.
\end{align}

The linearized heuristic \eqref{eq:linearizedHeuristic} defines a softmax decision-making model with a linear objective function of the form \eqref{eq:kroneckerModel}. Thus we can apply our estimation algorithm to estimate the parameters $\bftheta$. Using  \eqref{eq:linParamDef} we can then use the estimate of $\bftheta$ to provide an estimate of the parameters $(\mu_0,\sigma_0^2,\nu)$.

\subsection{Example estimations}
\label{sec:UCLexamples}
\begin{figure}
\centering
\includegraphics[width=3in]{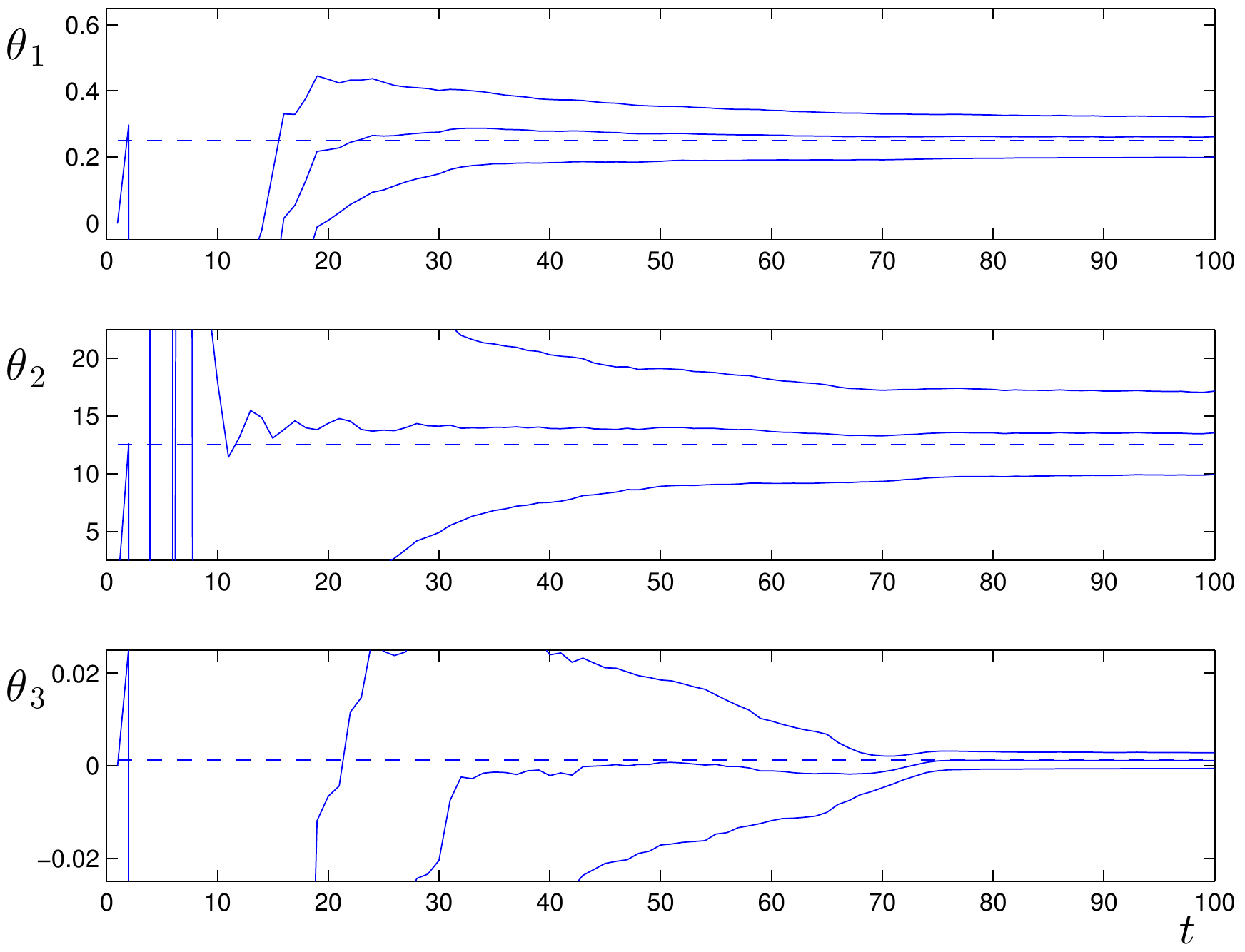}
\caption{Estimate of the vector of parameters $\bftheta$ based on simulated data from the stochastic UCL algorithm. The linearization point was taken to be $\bar{\mu}_0 = 150, \bar{\sigma}_0^2 = 2$. The true algorithm parameters were $\mu_0 = 200, \sigma_0^2 = 1, \lambda = 1$, and $\nu = 4$. The estimate converges as the number of observations $t$ grows. The dashed lines show the true value of each parameter $\theta_i$. For each value of $t$, an ensemble of 100 parameter estimates was formed by repeatedly simulating the data $\left\{(\bfx_t,\bfy_t)\right\}_{t=1}^T$ while holding the parameters $\bftheta$ fixed, and using the estimator to compute the value of the parameters. The solid lines show the mean parameter estimate and the 95\% confidence interval implied by the asymptotic normal distribution \eqref{eq:thetaMLAsymptotic}.}
\label{fig:fitNu}
\end{figure}

\begin{figure}
\centering
\includegraphics[width=3in]{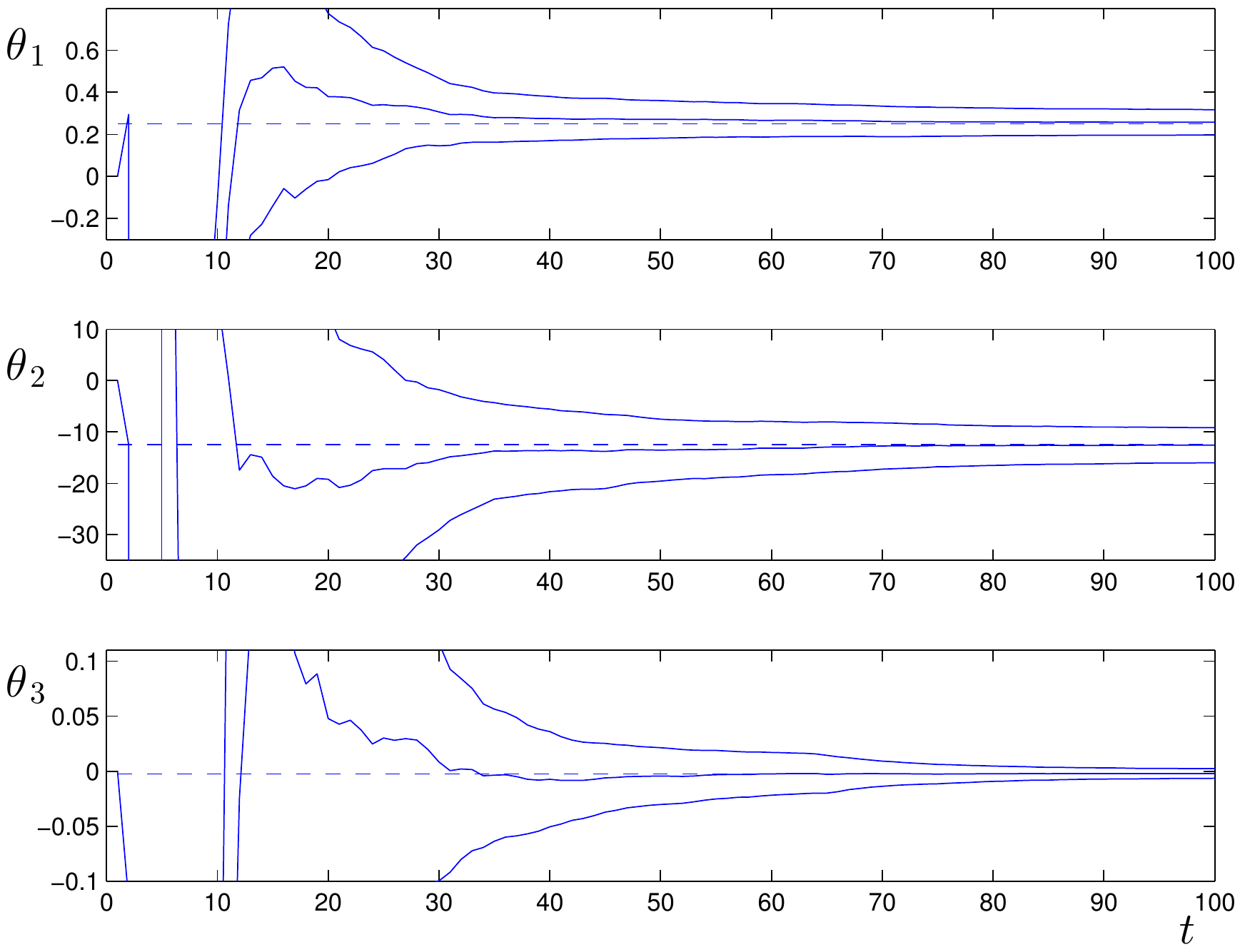}
\caption{Estimate of the vector of parameters $\bftheta$ based on simulated data from the stochastic UCL algorithm. Everything is the same as in Figure~\ref{fig:fitNu} except that the linearization point was taken to be $\bar{\mu}_0 = 250, \bar{\sigma}_0^2 = 0.5$. }
\label{fig:sim2b}
\end{figure}

We tested the estimation procedure described above by simulating runs of the stochastic UCL algorithm for various parameter values. Figures \ref{fig:fitNu} and \ref{fig:sim2b} show two examples of estimates computed using simulated data from the stochastic UCL algorithm with the nonlinear objective function $(\tilde{Q}_i^t \log t)/\nu$ and true parameters $(\mu_0,\sigma_0^2,\lambda,\nu) = (200,1,1,4)$. These parameters result in the algorithm achieving high performance (specifically, logarithmic regret, see \cite{PR-VS-NEL:14} for details). Figure \ref{fig:fitNu} shows estimates based on linearization about the point $(\muBar,\bar{\sigma}_0^2) = (150, 2)$. Following \eqref{eq:linParamDef}, the linearized objective function corresponds to parameters $\theta_1, \theta_2$, and $\theta_3$ having true values $\theta_1 = \frac{1}{\nu} = 0.25, \theta_2 = \frac{\mu_0-\bar{\mu}_0}{\nu} = 12.5$, and $\theta_3 = 1.25\times 10^{-3}$. These are the values to which the estimates should converge. Figure \ref{fig:sim2b} shows estimates based on linearization about the point $(\muBar,\bar{\sigma}_0^2) = (250,0.5)$. The linearized objective function in this case corresponds to the three parameters taking true values $\theta_1 = 0.25, \theta_2 = -12.5,$ and $\theta_3 = -2.5 \times 10^{-3}$.

In both cases the estimator converges to the true value of $\bftheta$ within the horizon $T=100$ of the decision task.  Further, the true value of the parameter is within the 95\% confidence interval after 30 observed choices. There are two implications from this result. First, the estimation procedure is at least somewhat robust to the choice of linearization point for this set of algorithm parameters. Second, the estimator is useful for realistic empirical data sets, such as those reported in \cite{PR-VS-NEL:14} and studied in the following section. For these data sets the horizon is $T=90$ choices. For this amount of data, the simulations show that the estimation procedure can identify the true value of the parameter in a statistically significant way. This result is valuable because the rigorous convergence result from Theorem \ref{thm:mlConvergence} does not directly guarantee convergence in the more general case of nonlinear objective functions.

\begin{figure}
\centering
\includegraphics[width=3in]{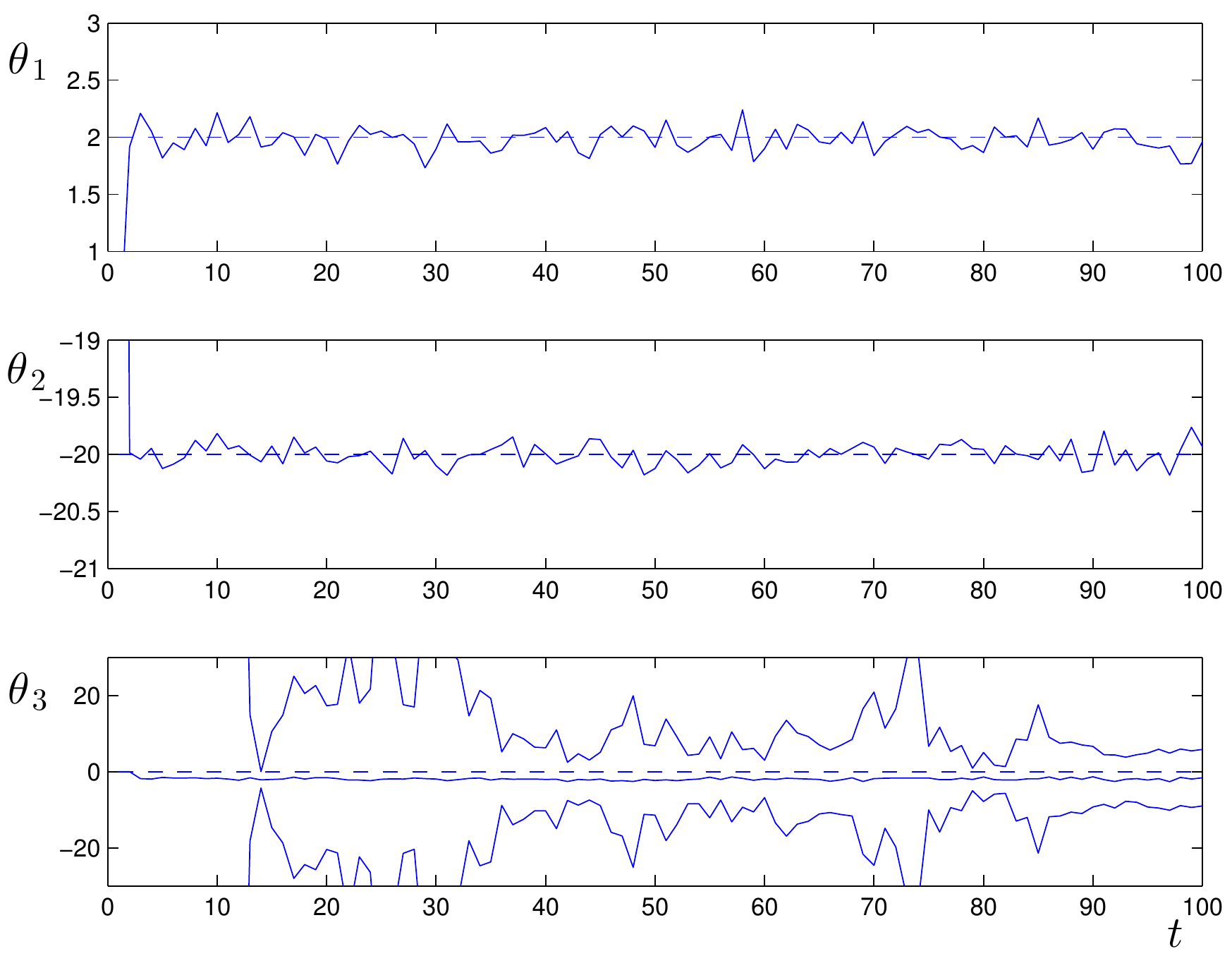}
\caption[Estimate of the vector of parameters $\bftheta$ based on simulated data from the stochastic UCL algorithm with a weakly-informative prior.]{Estimate of the vector of parameters $\bftheta$ based on simulated data from the UCL algorithm with a weakly-informative prior. This prior makes the algorithm's choice behavior more random, which makes the estimation problem more difficult. Everything is the same as in Figure~\ref{fig:fitNu} except that the linearization point was taken to be $\bar{\mu}_0 = 40, \bar{\sigma}_0^2 = 950$ and the true algorithm parameters were $\mu_0 = 30, \sigma_0^2 = 10^3, \lambda = 0$, and $\nu = 0.5$. 
The 95\% confidence interval implied by the asymptotic normal distribution \eqref{eq:thetaMLAsymptotic} is shown only in the plot of $\theta_3$.
For parameters $\theta_1$ and $\theta_2$, the width of the confidence intervals are much greater than the magnitudes of the parameter estimates and are omitted for legibility.}
\label{fig:fitLinRegret}
\end{figure}

The amount of data required to get a reliable estimate can depend on the true value of the algorithm parameters, as shown in Figure \ref{fig:fitLinRegret}. In this case, the true value of the algorithm parameters are $(\mu_0,\sigma_0^2,\lambda,\nu) = (30,10^3,0,0.5)$ and the linearization is made about the point $(\muBar,\bar{\sigma}_0^2) = (40,950)$. The linearized objective function   corresponds to the three parameters taking true values $\theta_1 = 2, \theta_2 = -20,$ and $\theta_3 = -1.05 \times 10^{-6}$. With the true values of the prior in the algorithm, the agent is sufficiently uncertain about the rewards and makes most of its initial 100 choices at random in order to gain information about the rewards. This choice behavior results in low performance (specifically, linear regret, see \cite{PR-VS-NEL:14} for details). Since the initial choices are effectively made at random, they do not provide useful information about the parameter values (except that they represent some combination of an uncertain prior and high decision noise). The uncertainty in the parameter values can be seen from the width of the confidence interval around the mean parameter estimates shown in Figure \ref{fig:fitLinRegret}. For $\theta_1$ and $\theta_2$ their width is many orders of magnitude larger than the magnitude of the parameter and they are not displayed. For $\theta_3$, the estimate exhibits persistent bias away from the true value, but the width of the associated confidence interval is significantly larger than the bias. Therefore, for such parameter values, one must observe more data to be able to shrink the confidence intervals and provide precise estimates of the parameter values.

\subsection{Discussion}
The linearization procedure described above yields a local linear approximation to the likelihood maximization problem \eqref{eq:mlProblem}, and Theorem \ref{thm:mlConvergence} provides conditions under which the local approximation results in an identified model with a convex optimization problem. However, the effectiveness of the procedure is sensitive to the choice of nominal prior $(\muBar,\deltaBar)$ about which to linearize.  The linearization point should be chosen such that the linear approximation is valid at the (unknown) true value of the parameters. In the worst case, there might not be any intuition for choosing the linearization point, making the above procedure no better than any other local optimization technique for which a starting point must be chosen.

Fortunately, there are several advantageous aspects of the problem. The first is generic to any heuristic function, which is the fact that the likelihood function forms a unique objective for judging the ``goodness'' of the estimated parameter. Without knowing in advance a good choice of linearization point, one approach is to perform the estimation assuming two different choices of linearization points and to compare the resulting estimates $\hat{\bftheta}$. If the two linearization points result in identical estimates there is no conflict, while if the estimates differ, the one with the higher likelihood value is better.

Second, there may be intuition about a appropriate choice of linearization point due to the structure of the model. In \cite{PR-VS-NEL:14}, we showed that behavior of the stochastic UCL model falls broadly into three classes as a function of the parameters $(\mu_0,\sigma_0^2,\lambda,\nu)$. Thus, by categorizing a given data set into one of the three classes, we narrow the search for a linearization point to the associated regions of parameter space. 
And, as we saw in Figures \ref{fig:fitNu} and \ref{fig:sim2b}, the stochastic UCL model appears to be relatively insensitive to the choice of linearization point within the region of parameter space associated with a given behavioral class. In the following section we exploit this intuition to estimate the parameters of the stochastic UCL algorithm based on data from a human subject experiment.

\section{Application to experimental data} \label{sec:experimentalData}
In this section we apply the estimator to fit the stochastic UCL model \eqref{eq:softmaxUCLPr} to experimental data studied in \cite{PR-VS-NEL:14}. By \emph{fit}, we refer to the process of selecting a nominal parameter for linearization and applying the estimator to the linearized model. The parameter estimates produced by the fitting procedure show that individuals with high performance match their behavior to the task in a statistically-significant way.

\subsection{Experimental setup}
This section reviews the experimental setup as presented in Reverdy \etal~\cite{PR-VS-NEL:14}. As described in \cite{PR-VS-NEL:14}, we collected data from a human subject experiment where we ran multi-armed bandit tasks through web servers at Princeton University (Princeton, NJ, USA) following protocols approved by the Princeton University Institutional Review Board. Human participants were recruited using Amazon's Mechanical Turk (AMT) web-based task platform \cite{MB-TK-SDG:11}. Participants were shown instructions that told them they would be playing a simple game during which they could obtain points, and that their goal was to obtain the maximum number of total points in each part of the game.

\begin{figure}
\centering
\includegraphics[width=3.0in]{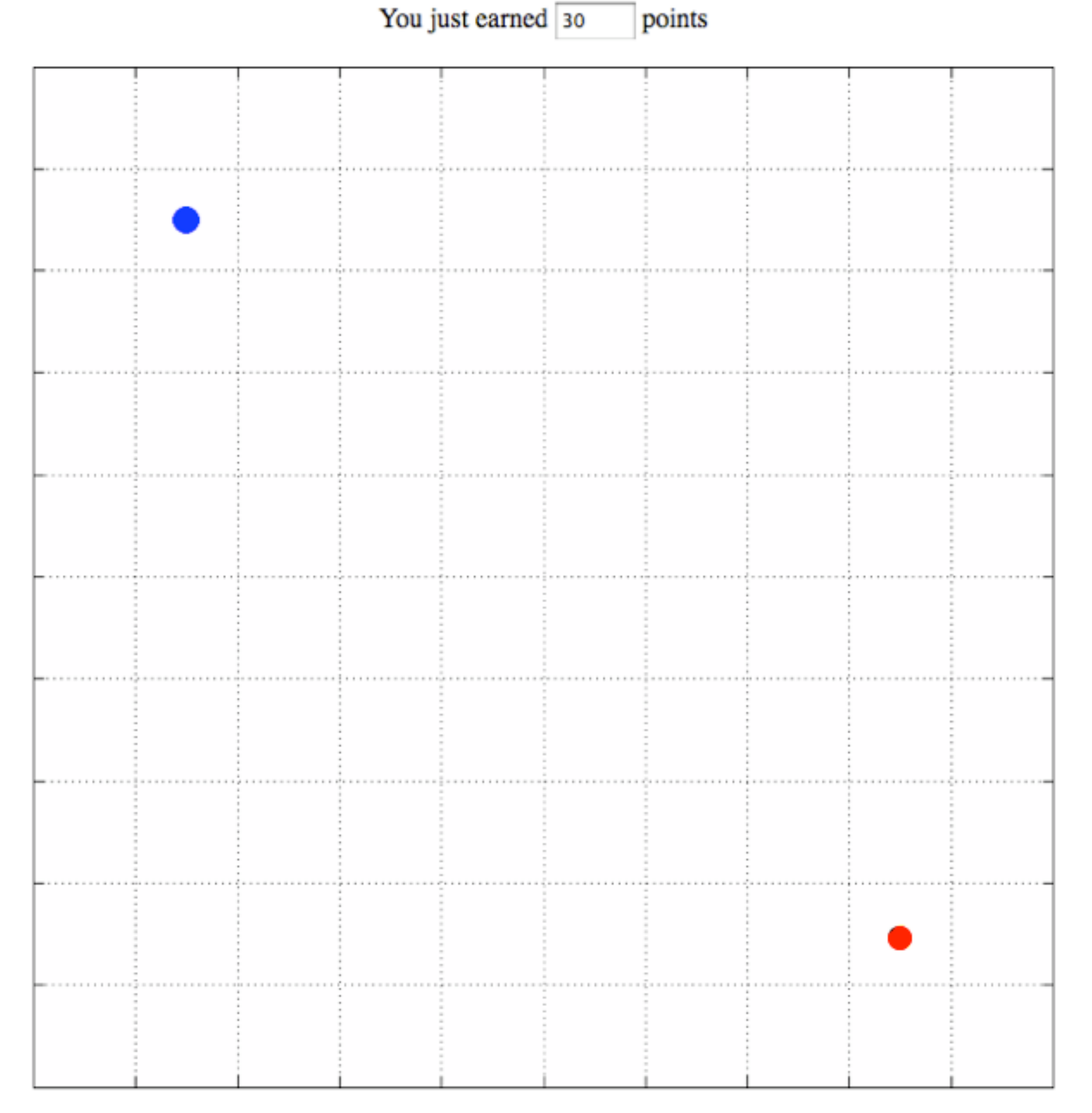}
\caption{The experimental interface used in the human subject experiment. Upon clicking on one of the 100 squares arranged in a $10 \times 10$ grid, the red dot would move to the center of the square. The subject was free to select a new square without penalty until the time allotted (1.5 or 6 seconds per choice) had elapsed, at which time the blue dot would move to the center of the selected square and the subject would receive a reward reported in the text box at the top of the screen. Originally appeared as Figure 5 of \cite{PR-VS-NEL:14}; reproduced with permission.
}
\vspace{-5mm}
\label{fig:surfaces}
\end{figure}

Each participant was presented with a set of $N=100$ options, presented as squares arranged in a $10 \times 10$ grid. See Figure \ref{fig:surfaces} for a visualization of the experimental interface. At each decision time $t \in \{1, \ldots, T\}$, the participant made a choice by moving the cursor to one square in the grid and clicking. After each choice was made, a numerical reward associated with that choice was reported on the screen. A variety of aspects of the game, including timing, game dynamics, and reward structures, were manipulated as part of the experimental design. As a result of these manipulations, only 326 of the 417 participants were assigned to a standard multi-armed bandit task for which the stochastic UCL model is appropriate. In the remainder of the section, we focus exclusively on data from these 326 participants.

\begin{figure}[h]
   \centering
   \begin{tabular}{c}
   	{\includegraphics[width=3.0in]{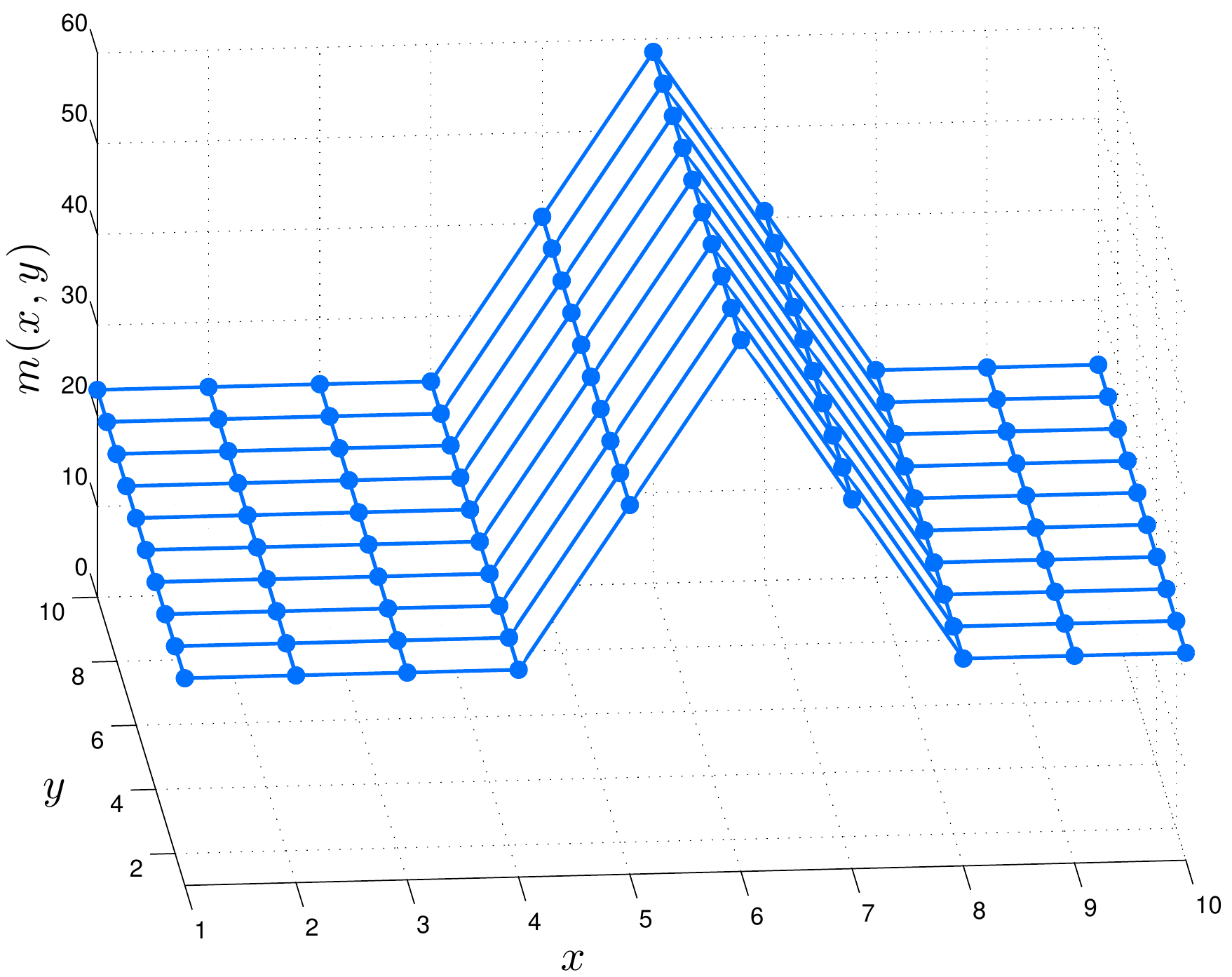}}\\
	(a)\\
   	{\includegraphics[width=3.0in]{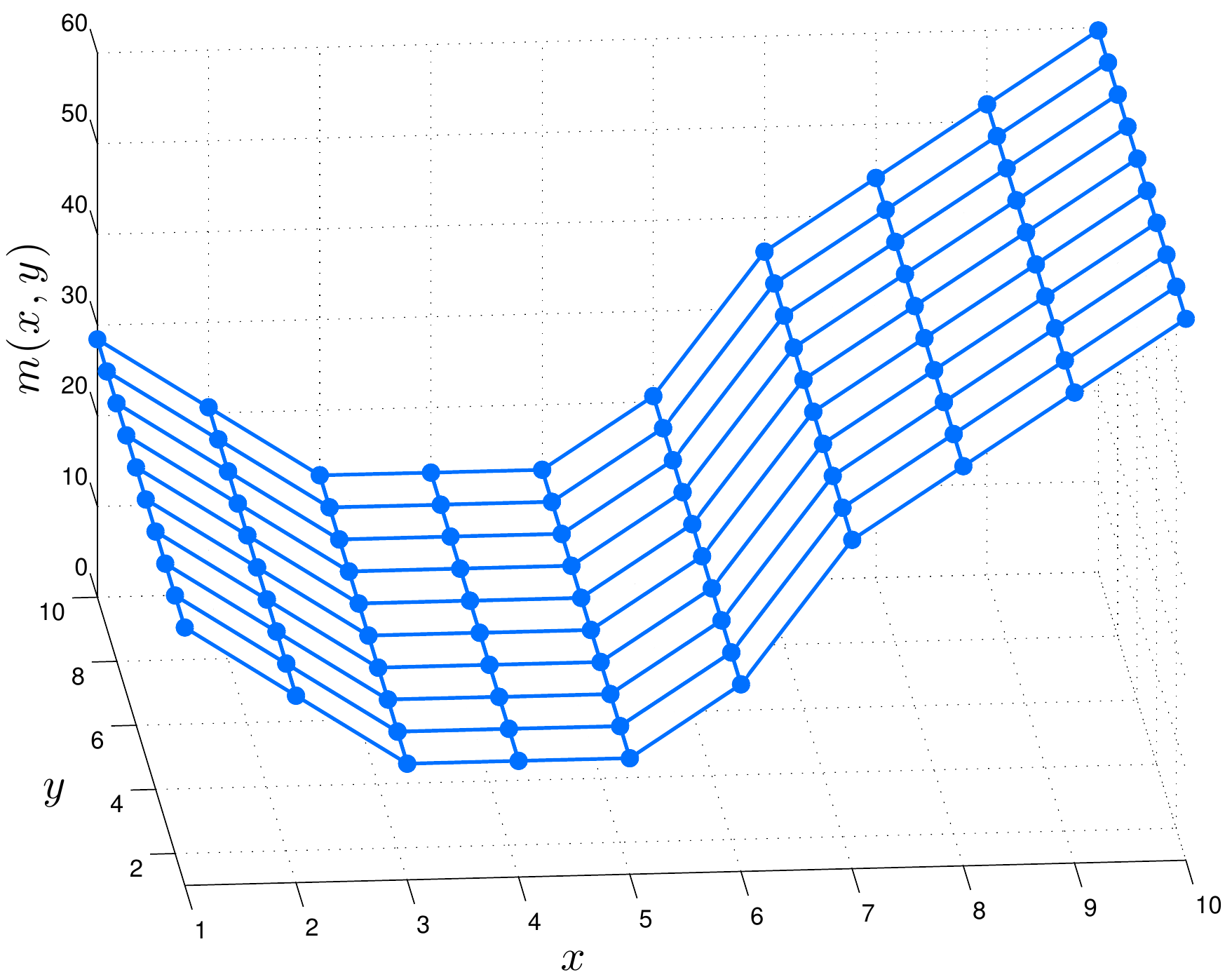}}\\
	(b)
  \end{tabular}
   \caption{The two task reward landscapes: (a) Landscape A, (b) Landscape B. The two-dimensional reward surfaces for the 10$\times$10 set of options followed the profile along one dimension (here the $x$ direction) and were flat along the other (here the $y$ direction). The Landscape A profile is designed to be simple in the sense that the surface is concave and there is only one global maximum ($x=6$), while the Landscape B profile is more complicated since it features two local maxima ($x = 1$ and 10), only one of which ($x=10$) is the global maximum. Originally appeared as Figure 6 of \cite{PR-VS-NEL:14}; reproduced with permission.}
\vspace{-6mm}
   \label{fig:surfaceProfiles}
\end{figure}

The mean value of the reward associated with choosing a particular option $i$ was $m_i$. Since the options were arranged in a $10 \times 10$ grid, the set of mean values can be thought of as a real-valued function on the discrete two-dimensional grid. We refer to this function as the reward landscape, and prior knowledge about the rewards in a given task corresponds to prior knowledge about the landscape. Mean rewards in each task corresponded to one of two landscapes: Landscape A and Landscape B, shown in Figure \ref{fig:surfaceProfiles}. Each landscape was flat along one dimension and followed a profile along the other dimension. The profile of Landscape A was such that a simple gradient-climbing strategy was likely to prove effective, while Landscape B was constructed to require a more sophisticated strategy. Each participant played the game with each landscape once, presented in random order.  Due to the structure of the experimental design, only one of the two landscapes was associated with a standard multi-armed bandit task.

The participants' performance in a given task can be classified in terms of the growth rate of their cumulative regret, which is a measure of cumulative loss relative to the (unknown to the subject) optimal decision. As reported in \cite{PR-VS-NEL:14}, 70 of the 326 participants, or approximately 21\%, achieved high performance while the remainder, approximately 79\%, achieved low performance. Of the 206 subjects assigned to Landscape A, 53 achieved high performance.  Likewise, of the 120 subjects assigned to Landscape B, 17 achieved high performance. The high-performing subjects outperformed standard frequentist algorithms on the task, which we attribute to the subjects' having good priors about the task. Since we did not explicitly convey prior knowledge about the reward landscapes to the subjects, we postulate that they used priors developed in the course of other spatial tasks encountered in daily life. Considering the stochastic UCL algorithm as a model of the subjects' behavior, good priors correspond to good values for the parameters $\bfmu_0$ and $\Sigma_0$, which quantify the subjects' intuition about the task.  To learn the priors we propose estimating them from the data. The estimated priors could then be used, e.g., to improve the performance of an automated system.

\subsection{Fitting}
In fitting the stochastic UCL model to human subject data, we seek to answer two questions. First, what distinguishes the decision-making of the subjects with high performance from those with low performance? And second, do subjects adapt their decision-making strategies to the task, i.e., the reward landscape? Our experimental design provides data from only one task per subject, so we cannot, for example, compare a single subject's performance on the different landscapes. Thus, we analyze at the  population-level   to answer the two questions.

Each subject is classified as having high or low performance as described above. On the basis of this classification and the reward landscape, the subject is assigned to one of the four performance--landscape combined categories. We assume each subject represents an independent and identically distributed (iid) sample from the true parameter $\bftheta_0$ associated with its category. We applied the estimator to data from each subject using nominal parameters $(\bar{\mu}_0,\bar{\sigma}_0^2,\lambda) = (30,10,0.1)$ for subjects with low performance and $(\bar{\mu}_0,\bar{\sigma}_0^2,\lambda) = (200,10^6,4)$ for subjects with high performance. We validated the choice of $\lambda$ by performing estimation on the data from several subjects using a variety of values of $\lambda$. The optimal value of $\lambda$ clearly differed between the two categories of performance but the estimates for each given performance category were fairly robust to changes in the value of $\lambda$. The fitting procedure produces a maximum likelihood estimate and associated covariance matrix for each subject. By the iid assumption, it is tenable to construct a population-level parameter estimate for each of the four categories by appropriately averaging the individual subjects' estimates and covariances.

Table \ref{tab:thetaHat} presents the population-level parameter estimates, along with the mean log likelihood values, for the four categories. The columns labeled $\hat{\theta}$ report the maximum likelihood parameter estimates and those labeled $\sigma$ their asymptotic standard deviations implied by \eqref{eq:thetaMLAsymptotic}. Recall that these parameters represent deviations from the nominal parameter values and therefore are not directly comparable between performance categories. However,  comparing the magnitude of the standard deviations shows that the parameter estimates are much more precise for those categories associated with high performance.   This is consistent with our findings in Section~\ref{sec:UCLexamples}. Table \ref{tab:thetaHatTransformed} presents the maximum likelihood parameter estimates transformed back into the original variables $\nu, \mu_0$, and $\sigma_0^2$; these are directly comparable.

Table \ref{tab:thetaHatTransformed} allows us to answer our first question about the differences between subjects with different levels of performance. The parameter values clearly differ more between levels of performance rather than between landscapes. Between levels of performance the parameters that differ the most are the decision noise parameter $\nu$ and the prior uncertainty $\sigma_0^2$. Larger values of $\nu$ are associated with more random decision-making, while larger values of $\sigma_0^2$ represent greater uncertainty about the rewards which is associated with placing a higher value on information. Both of these factors tend to encourage exploration, and the values of both $\nu$ and $\sigma_0^2$ are much greater for subjects with high performance than those with low performance. Thus, for both landscapes, the high-performing subjects explore more than the low-performing ones, which presumably helps them discover the regions of high rewards. Furthermore, the subjects with high performance use correlated priors which allow them to quickly explore large regions of the reward surface. 

We can compare the quality of the model fits by comparing the mean log likelihood values across categories provided on Table \ref{tab:thetaHat}. Again, we see starker differences between levels of performance than between landscapes.  Between landscapes, the fits are approximately equal in quality, while between performance levels there is substantial difference, equal to an approximate doubling of the fitted model's predictive power.

Table \ref{tab:thetaHat} allows us to answer our second question about the degree to which subjects match their strategies to the task. We focus on comparing the parameters across landscape conditions for each of the performance categories separately. For low-performing subjects, comparing the relative magnitudes of the parameter estimates and their standard deviations suggests that there is no significant difference between the two landscape conditions. The two-sided Welch's $t$-test \cite{BLW:47} confirms that the difference in the parameter estimates is statistically insignificant. For high-performing subjects, the parameter estimates are much more precise, and the two-sided Welch's $t$-test confirms that the difference in the parameter estimates is statistically significant at the 95\% confidence level. In other words, the fitting procedure is able to distinguish that the high-performing subjects have strategies that are matched to the landscape.

\subsection{Implications for human-centered automation}
The results of the fitting exercise have several implications for human-centered automation. First, they demonstrate an estimator for a model of human decision-making behavior. The estimator allows one to quantify a human subject's intuition in a statistically powerful way. Second, the model fits are of higher quality for subjects with high performance. This suggests that the stochastic UCL model is better suited to the decision-making behavior of subjects who are experts at the task; a different model may be more appropriate for lower-performing subjects. Third, subjects with high performance seem to have effective priors: these priors have low certainty (large values of $\sigma_0^2$), but exploit correlation in the rewards due to the smoothness of the reward landscapes by using positive values of the length scale parameter $\lambda$. When such correlation structures exist, they can be exploited to greatly improve performance \cite{VS-PR-NEL:15}, as our human subjects appear to have done. The estimator provides a way to learn effective priors from a human operator. In the absence of a correlation structure, the above fitting process can still be applied by setting $\lambda = 0$, although convergence of the estimator will be slower, requiring longer series of choice data than those studied here.

By analyzing data from a human subject experiment, we have shown the effectiveness of the linearization procedure for extending the estimator to a model with a nonlinear objective function. The known asymptotic properties of the estimator allowed us to perform tests for statistical significance and find differences in behavior.

\begin{table}
\centering
\begin{tabular}{cc|cc}
\multicolumn{4}{c}{Low (linear, power-law) performance} \\
\multicolumn{2}{c}{Landscape A, $153$ subj.} & \multicolumn{2}{c}{Landscape B, $103$ subj.}\\
\multicolumn{2}{c}{Mean log likelihood: -338} & \multicolumn{2}{c}{Mean log likelihood: -331}\\
$\hat{\theta}$ 	& $\sigma$ & $\hat{\theta}$ & $\sigma$ \\
0.360	& 90.4 	& 0.252	& 1.32 \\
-5.22 	& 1.27e3 	& -2.12	& 51.8 \\
0.433 	& 1.02e2 	& 0.213 	& 8.61 \\
\end{tabular}

\begin{tabular}{cc|cc}
\multicolumn{4}{c}{High (log-law) performance} \\
\multicolumn{2}{c}{Landscape A, $53$ subj.} & \multicolumn{2}{c}{Landscape B, $17$ subj.}\\
\multicolumn{2}{c}{Mean log likelihood: -273} & \multicolumn{2}{c}{Mean log likelihood: -271}\\
$\hat{\theta}$ & $\sigma$ & $\hat{\theta}$ & $\sigma$ \\
3.93e-2 	& 1.18e-3 	& 3.39e-2 	& 1.04e-3 \\
-6.86 	& 0.226 	& -6.57	& 0.268 	\\
7.88e-7 	& 2.34e-8 	& 6.80e-7 	& 2.06e-8	\\
\end{tabular}
\caption{Parameter estimates $\hat{\theta}$ and associated standard deviations $\sigma$ conditional on regret growth order and reward landscape. The values for high performance are significantly different between surfaces at the 95\% confidence level (two-sided Welch's $t$-test \cite{BLW:47}); other comparisons show that the parameter values do not significantly differ between classes.}
\label{tab:thetaHat}
\end{table}

\begin{table}
\centering
\begin{tabular}{cc|cc}
\multicolumn{4}{c}{Low (linear, power-law) performance} \\
\multicolumn{2}{c}{Landscape A, $153$ subj.} & \multicolumn{2}{c}{Landscape B, $103$ subj.}\\
Parameter 	& Value	& Parameter 	& Value \\
$\nu$ 		& 2.78  	& $\nu$ 		& 3.97 \\
$\mu_0$ 		& 15.5  	& $\mu_0$ 	& 21.6 \\
$\sigma_0^2$	& 4.54 	& $\sigma_0^2$ & 5.42 \\
\end{tabular}

\begin{tabular}{cc|cc}
\multicolumn{4}{c}{High (log-law) performance} \\
\multicolumn{2}{c}{Landcape A, $53$ subj.} & \multicolumn{2}{c}{Landscape B, $17$ subj.}\\
Parameter 	& Value	& Parameter 	& Value \\
$\nu$ 		& 25.5  	& $\nu$ 		& 29.5 \\
$\mu_0$ 		& 25.3 	& $\mu_0$ 	& 6.08 \\
$\sigma_0^2$	& 3.32e5 	& $\sigma_0^2$ & 3.35e5 \\
\end{tabular}
\caption{Parameter estimates $\nu, \mu_0, \sigma_0^2$ and associated standard deviations $\sigma$ conditional on regret growth order and reward landscape.}
\label{tab:thetaHatTransformed}
\vspace{-10mm}
\end{table}

\section{Conclusion}
Motivated by the parameter estimation problem for decision-making models, we studied the maximum likelihood parameter estimation problem for softmax decision-making models with linear objective functions. Such models occur frequently in the neuroscience and machine learning literatures. We derived conditions under which the maximum likelihood estimator converges on the correct parameter values, characterized the estimator's asymptotic distribution, and showed how to use this distribution to formulate confidence intervals for the parameter estimates.

We then showed that the stochastic UCL algorithm could be transformed into a softmax decision-making model with a linear objective function by linearizing the objective function about a nominal point in parameter space. By performing parameter estimation on the linearized model using simulated data, we showed that we could estimate the true value of the stochastic UCL algorithm parameters. The amount of data required to perform useful estimation depends on the region of parameter space, with parameters representing priors that strongly influence behavior (e.g., small variances $\sigma_0^2$, representing strong beliefs, or large correlation length scales $\lambda$, representing highly structured beliefs) being easier to estimate. This is intuitive, as observed behavior will be more sensitive to such influential beliefs.

The estimator convergence results we state in Theorem \ref{thm:mlConvergence} 
are specific to the case where the objective function is a linear function of the unknown parameters. However, we showed how the estimation procedure can be extended to nonlinear objective functions with linearization. Using the linearization technique with the estimator, we fit the stochastic UCL model developed in \cite{PR-VS-NEL:14} to data from a human subject experiment. The estimates show a statistically significant difference in behavior between subjects who exhibit good performance in similar but different tasks. Quantifying these differences are of interest both for the science of decision-making but also for the development of automation technology. In conjunction with the stochastic UCL model, the estimator developed in this paper provides the tools for quantifying human decision-making behavior in multi-armed bandit problems. These tools will facilitate the principled development of human-machine decision-making teams.

\bibliography{glm-revised}

\begin{thebibliography}{10}

\bibitem{DB:92}
D.~B{\"o}hning.
\newblock Multinomial logistic regression algorithm.
\newblock {\em Annals of the Institute of Statistical Mathematics},
  44(1):197--200, 1992.

\bibitem{CGB:70}
C.~G. Broyden.
\newblock The convergence of a class of double-rank minimization algorithms.
\newblock {\em IMA Journal of Applied Mathematics}, 6(1):76--90, 1970.

\bibitem{MB-TK-SDG:11}
M.~Buhrmester, T.~Kwang, and S.~D. Gosling.
\newblock {A}mazon's {M}echanical {T}urk: A new source of inexpensive, yet
  high-quality, data?
\newblock {\em Perspectives on Psychological Science}, 6(1):3--5, 2011.

\bibitem{JDC-SMM-AJY:07}
J.~D. Cohen, S.~M. McClure, and J.~Y. Angela.
\newblock Should {I} stay or should {I} go? how the human brain manages the
  trade-off between exploitation and exploration.
\newblock {\em Philosph. Trans. of the Roy. Soc. B: Biological Sci.},
  362(1481):933--942, 2007.

\bibitem{NDD:11}
N.~D. Daw.
\newblock Trial-by-trial data analysis using computational models.
\newblock {\em Decision making, affect, and learning: Attention and performance
  XXIII}, 23:3--38, 2011.

\bibitem{NDD-etal:06}
N.~D. Daw, J.~P. O'Docherty, P.~Dayan, B.~Seymour, and R.~J. Dolan.
\newblock Cortical substrates for exploratory decisions in humans.
\newblock {\em Nature}, 441(7095):876--879, 2006.

\bibitem{RF:70}
R.~Fletcher.
\newblock A new approach to variable metric algorithms.
\newblock {\em The {C}omput. {J}.}, 13(3):317--322, 1970.

\bibitem{KG-NAS:10a}
K.~Gimpel and N.~A. Smith.
\newblock Softmax-margin training for structured log-linear models.
\newblock Technical Report CMU-LTI-10-008, Carnegie Mellon Univ., 2010.

\bibitem{JG-etal:10}
J.~Gl{\"a}scher, N.~Daw, P.~Dayan, and J.~P. O'Doherty.
\newblock States versus rewards: dissociable neural prediction error signals
  underlying model-based and model-free reinforcement learning.
\newblock {\em Neuron}, 66(4):585--595, 2010.

\bibitem{ASG:91}
A.~S. Goldberger.
\newblock {\em A Course in Econometrics}.
\newblock Harvard Univ. Press, 1991.

\bibitem{DG:70}
D.~Goldfarb.
\newblock A family of variable-metric methods derived by variational means.
\newblock {\em Mathematics of {C}omputation}, 24(109):23--26, 1970.

\bibitem{SMK:93}
S.~Kay.
\newblock {\em Fundamentals of Statistical Signal Processing, Volume I:
  Estimation Theory}.
\newblock Prentice Hall, 1993.

\bibitem{BK-LC-MATF-AJH:05}
B.~Krishnapuram, L.~Carin, M.~A.~T. Figueiredo, and A.~J. Hartemink.
\newblock Sparse multinomial logistic regression: Fast algorithms and
  generalization bounds.
\newblock {\em IEEE Trans. Pattern Anal. Mach. Intell.}, 27(6):957--968, 2005.

\bibitem{BL-PWG:05}
B.~Lau and P.~W. Glimcher.
\newblock Dynamic response-by-response models of matching behavior in rhesus
  monkeys.
\newblock {\em Journal of the Experimental Analysis of Behavior},
  84(3):555--579, Nov 2005.

\bibitem{DM:74}
D.~McFadden.
\newblock Conditional logit analysis of qualitative choice behavior.
\newblock In P.~Zarembka, editor, {\em Frontiers in Econometrics}, pages
  105--142. Academic Press, New York, 1974.

\bibitem{DM-FR-ASV:86}
D.~Mitra, F.~Romeo, and A.~Sangiovanni-Vincentelli.
\newblock Convergence and finite-time behavior of simulated annealing.
\newblock {\em Advances in Applied Probability}, 18(3):747--771, 1986.

\bibitem{PRM-BKC-JDC:06}
P.~R. Montague, B.~King-Casas, and J.~D. Cohen.
\newblock Imaging valuation models in human choice.
\newblock {\em Annu. Rev. Neurosci.}, 29:417--448, 2006.

\bibitem{MRN-JIG:13}
M.~R. Nassar and J.~I. Gold.
\newblock A healthy fear of the unknown: Perspectives on the interpretation of
  parameter fits from computational models in neuroscience.
\newblock {\em PLoS Computational Biology}, 9(4):e1003015, 2013.

\bibitem{AN-etal:12}
A.~Nedic, D.~Tomlin, P.~Holmes, D.~A. Prentice, and J.~D. Cohen.
\newblock A decision task in a social context: Human experiments, models, and
  analyses of behavioral data.
\newblock {\em Proc. IEEE}, 100(3):713--733, 2012.

\bibitem{WKN-DLM:94}
W.~K. Newey and D.~McFadden.
\newblock Large sample estimation and hypothesis testing.
\newblock In R.~F. Engle and D.~L. McFadden, editors, {\em Handbook of
  Econometrics}, volume~4, chapter~36, pages 2111--2245. Elsevier, 1994.

\bibitem{AYN-SJR:00}
A.~Y. Ng and S.~J. Russell.
\newblock Algorithms for inverse reinforcement learning.
\newblock In {\em The Int. Conf. on Machine Learning}, pages 663--670, 2000.

\bibitem{VR-HB:11}
V.~Ramanujam and H.~Balakrishnan.
\newblock Estimation of maximum-likelihood discrete-choice models of the runway
  configuration selection process.
\newblock In {\em Proc. American Control Conf.}, pages 2160--2167, 2011.

\bibitem{ReverdyThesis}
P.~Reverdy.
\newblock {\em Human-inspired algorithms for search: A framework for
  human-machine multi-armed bandit problems}.
\newblock PhD thesis, Princeton University, Department of Mechanical and
  Aerospace Engineering, 2014.

\bibitem{PR-VS-NEL:14}
P.~Reverdy, V.~Srivastava, and N.~E. Leonard.
\newblock Modeling human decision-making in generalized {G}aussian multi-armed
  bandits.
\newblock {\em Proc. IEEE}, 102(4):544--571, 2014.

\bibitem{HR:52}
H.~Robbins.
\newblock Some aspects of the sequential design of experiments.
\newblock {\em Bulletin of the Amer. Math. Soc.}, 58:527--535, 1952.

\bibitem{SR:98}
S.~Russell.
\newblock Learning agents for uncertain environments.
\newblock In {\em Proc. 11th Annu. Conf. on Computational Learning Theory},
  pages 101--103. ACM, 1998.

\bibitem{KS-etal:05}
K.~Samejima, Y.~Ueda, K.~Doya, and M.~Kimura.
\newblock Representation of action-specific reward values in the striatum.
\newblock {\em Science}, 310(5752):1337--1340, 2005.

\bibitem{DFS:70}
D.~F. Shanno.
\newblock Conditioning of quasi-{N}ewton methods for function minimization.
\newblock {\em Mathematics of {C}omputation}, 24(111):647--656, 1970.

\bibitem{VS-PR-NEL:15}
V.~Srivastava, P.~Reverdy, and N.~E. Leonard.
\newblock Correlated multiarmed bandit problem: Bayesian algorithms and regret
  analysis.
\newblock {\em arXiv:1507.01160v2}, 2015.
\newblock Submitted for publication.

\bibitem{ARS-etal:12}
A.~R. Stewart, M.~Cao, A.~Nedic, D.~Tomlin, and N.~E. Leonard.
\newblock Towards human-robot teams: Model-based analysis of human decision
  making in two-alternative choice tasks with social feedback.
\newblock {\em Proc. IEEE}, 100(3):751--775, 2012.

\bibitem{RSS-AGB:98}
R.~S. Sutton and A.~G. Barto.
\newblock {\em Introduction to reinforcement learning}.
\newblock MIT Press, 1998.

\bibitem{CW-PD:92}
C.~J. C.~H. Watkins and P.~Dayan.
\newblock ${Q}$-learning.
\newblock {\em Machine learning}, 8(3-4):279--292, 1992.

\bibitem{BLW:47}
B.~L. Welch.
\newblock The generalization of ``{S}tudent's'' problem when several different
  population variances are involved.
\newblock {\em Biometrika}, 34(1-2):28--35, 1947.

\bibitem{RCW-YN:13}
R.~C. Wilson and Y.~Niv.
\newblock Is model fitting necessary for model-based f{MRI}?
\newblock In {\em Proceedings of the Multi-disciplinary Conference on
  Reinforcement Learning and Decision Making}, page S41, 2013.

\end{thebibliography}
\bibliographystyle{abbrv}

\end{document}